\documentclass[10pt,letterpaper]{amsart}
\usepackage{iftex}
\ifPDFTeX
  \usepackage[T1]{fontenc}
  \usepackage[utf8]{inputenc}
  \usepackage{lmodern}
\else
  \usepackage{fontspec}
\fi
\usepackage{amsmath,amssymb,mathrsfs,mathtools}
\usepackage[textwidth=444pt,top=0.9in,bottom=0.9in]{geometry}
\usepackage{enumitem}
\usepackage{xcolor}
\usepackage[unicode,colorlinks=true,linkcolor=green!35!black,citecolor=green!35!black,urlcolor=green!35!black]{hyperref}
\usepackage{bookmark}
\setlist[enumerate]{label=\textup{(\arabic*)},leftmargin=2.5em,itemsep=2pt,topsep=3pt}
\numberwithin{equation}{section}
\newtheorem{thm}{Thm}[section]
\newtheorem{prop}[thm]{Prop}
\newtheorem{cor}[thm]{Cor}
\theoremstyle{definition}
\newtheorem{defn}{Def}[section]
\newtheorem{exmp}[thm]{Example}
\newtheorem{rem}[thm]{Remark}
\newtheorem{ques}[thm]{Question}
\newcommand{\Z}{\mathbb Z}
\newcommand{\Q}{\mathbb Q}
\newcommand{\C}{\mathbb C}
\newcommand{\N}{\mathbb N}
\newcommand{\Gr}{\mathrm{Gr}}
\newcommand{\Sh}{\mathrm{Sh}}
\newcommand{\Ab}{\mathrm{Ab}}
\newcommand{\Tor}{\operatorname{Tor}}
\newcommand{\Ext}{\operatorname{Ext}}
\newcommand{\Hom}{\operatorname{Hom}}
\newcommand{\Ann}{\operatorname{Ann}}
\newcommand{\fd}{\operatorname{fd}}
\newcommand{\supp}{\operatorname{supp}}
\newcommand{\id}{\mathrm{id}}
\newcommand{\ord}{\mathcal O}
\newcommand{\sh}[1]{\langle #1\rangle}
\newcommand{\HH}{\mathcal H}
\newcommand{\cN}{\mathcal N}
\newcommand{\cE}{\mathcal E}
\newcommand{\eps}{\varepsilon}
\newcommand{\up}{\uparrow}

\newcommand{\byref}[2]{\href{https://arxiv.org/pdf/2610.38300v1\#page=#1}{[BY, #2]}}
\title[Homology and Topology of Partially Ordered Fields]{Homology and Topology of Partially Ordered Fields}
\author{Hanzhen Bao}
\date{}
\hypersetup{pdftitle={Homology and Topology of Partially Ordered Fields},pdfauthor={Hanzhen Bao},pdfsubject={English manuscript}}
\makeatletter
\def\@secnumfont{\bfseries}
\def\@settitle{\begin{center}
  \normalfont\fontsize{16}{20}\selectfont\bfseries
  \@title\par
\end{center}}
\def\@setauthors{\begin{center}
  \normalfont\fontsize{11}{14}\selectfont\authors
\end{center}}
\def\part{\@startsection{part}{0}
  \z@{1.3\linespacing\@plus.5\linespacing}{.7\linespacing}
  {\normalfont\fontsize{14}{18}\selectfont\bfseries\boldmath\raggedright}}
\def\section{\@startsection{section}{1}
  \z@{1.1\linespacing\@plus.5\linespacing}{.6\linespacing}
  {\normalfont\fontsize{13}{16}\selectfont\bfseries\boldmath\centering}}
\def\specialsection{\@startsection{section}{1}
  \z@{1.1\linespacing\@plus.5\linespacing}{.6\linespacing}
  {\normalfont\fontsize{13}{16}\selectfont\bfseries\boldmath\centering}}
\def\subsection{\@startsection{subsection}{2}
  \z@{.9\linespacing\@plus.3\linespacing}{.4\linespacing}
  {\normalfont\fontsize{11.5}{14}\selectfont\bfseries\boldmath\raggedright}}
\def\subsubsection{\@startsection{subsubsection}{3}
  \z@{.7\linespacing\@plus.3\linespacing}{.35\linespacing}
  {\normalfont\fontsize{11}{14}\selectfont\bfseries\boldmath\raggedright}}

\renewenvironment{abstract}{%
  \global\setbox\abstractbox=\vtop\bgroup
  \normalfont\small
  \list{}{\leftmargin1.5em\rightmargin1.5em
    \listparindent\normalparindent\itemindent\z@\parsep\z@
    }
  \item\relax
  {\centering\normalfont\fontsize{11.5}{14}\selectfont\bfseries\abstractname\par}
  \vspace{.4\baselineskip}\noindent\ignorespaces
}{\endlist\egroup}
\makeatother

\begin{document}
\begin{abstract}
We develop a unified homology theory for directed partial orders (DPOs) and lattice-ordered positive cones on fields of characteristic zero. With the underlying field and its operations fixed, we define order-carrying spaces independently by topological conditions, establish a bijection with DPO positive cones, and characterize lattice orders by intersections of minimal open neighborhoods. A unified comparison theorem identifies graded $\Ext$ for region ideals with sheaf cohomology and multiplicative $\Tor$ with the reduced singular homology of truncated regions, giving a topological interpretation of the existence and preservation of reduced decompositions. We obtain canonical cohomological obstructions to lattice ordering, prove that a coefficient module for positive multiplication is flat if and only if the multiplier has a positive inverse, and show that every strict extension of positive cones on the same field produces a nonzero first $\Tor$ obstruction. Finally, we place the concrete lattice recently constructed by Bao--Zhang within this framework.
\end{abstract}
\maketitle
\tableofcontents
\bigskip
\noindent\textbf{Note on This Preliminary Draft}\par\smallskip
{\bfseries
This manuscript is a very rough first draft on the subject. It was written in its entirety by GPT6Astra under the author's guidance in natural language, to present the ideas as quickly as possible. The exposition of the theoretical framework, the statements, and the proofs remain far from fully refined. The manuscript will be revised and updated on an ongoing basis. Please contact the author with any errors or concerns.\par}
\medskip
\section*{Introduction}
\subsection*{Background}
In 1956, Birkhoff and Pierce asked whether the complex field admits a lattice order compatible with addition and multiplication \cite{BP}. Schwartz and Yang proved the existence of directed partial orders (DPOs) on the complex field in 2011 and subsequently gave necessary and sufficient conditions for a general field to admit a compatible directed partial order \cite{SY11,SY23}.

Recently, Wang, Yuan, Zhang, and Zhu constructed an inverse-closed directed positive cone on the complex field using integral closures of local rings \cite{WYZZ}. Subsequently, He and Bao--Zhang independently constructed compatible lattice orders on the complex field \cite{He,BY}. Bao--Zhang also characterized lattice ordering by the existence of a reduced decomposition for every element and used valuation-controlled transfinite extensions to construct and preserve these decompositions on the fixed complex field.

These advances bring structural questions into focus: when is a directed partial order a lattice order? How does positive multiplication affect lattice operations? Can reduced decompositions be preserved when a positive cone is enlarged? Which intrinsic properties distinguish lattice orders arising from different constructions? The preceding existence and construction results have not placed these questions within a unified homological framework. Starting from decompositions as differences of positive elements, we develop a language connecting algebra, topology, and homology in which these questions can be formulated, computed, and compared.

\subsection*{Main results}
We first establish a correspondence between positive cones and order-carrying spaces. On a fixed field $L$ of characteristic zero, we define an order-carrying space as a connected $T_0$ Alexandrov topology for which translations are homeomorphisms and field multiplication is continuous at the origin. This definition takes neither a positive cone nor a DPO as input. We prove that order-carrying topologies correspond bijectively to DPO positive cones, with the cone recovered as the minimal open neighborhood of zero. Lattice ordering is expressed by an intrinsic intersection property: the intersection of any two minimal open neighborhoods is again the minimal open neighborhood of a point. It can also be characterized by the doubling map $x\mapsto2x$ being a homeomorphism and by compactness of these intersections. Thus, with the underlying field and its operations retained, the algebraic axioms of a positive cone receive a complete spatial formulation, and the existence of suprema becomes recognizable from the topology.

The other central result is the unified comparison theorem. Let $P$ be a DPO positive cone on $L$, set $R=\Z[P]$, and denote its order-carrying space by $X_P$. We prove an exact equivalence between the category of $L$-graded $R$-modules and the category of sheaves of abelian groups on $X_P$, identifying graded $\Ext$ for region ideals with ordinary sheaf cohomology. For a positive multiplier $a\in P\setminus\{0\}$, let $R_a$ be the coefficient module whose left action is twisted by $X^p\mapsto X^{ap}$. For an upward-closed subset $V\subseteq P$ and $u\in P$, write
\[
J_V=\bigoplus_{v\in V}\Z X^v,\qquad
Y_{a,u}(V)=\{v\in V:av\le_Pu\}\subseteq X_P.
\]
There are natural isomorphisms
\[
\Tor_n^R(R/J_V,R_a)_u
\cong\widetilde H_{n-1}\!\left(Y_{a,u}(V);\Z\right),\qquad n\ge1.
\]
The right-hand side is the usual singular homology of the truncated region, and positive translations correspond to maps induced by inclusions of regions. Algebraic relations and their higher obstructions can therefore be tracked directly through connected components and higher homology.

This framework also characterizes the existence of reduced decompositions: the canonical cohomological obstruction of an element vanishes if and only if that element has a reduced decomposition. Simultaneous vanishing of all these obstructions characterizes lattice ordering.

Having characterized existence, we study preservation under positive multiplication and extension of positive cones. Although positive multiplication preserves order, it may introduce a new common positive part into a reduced pair. For a lattice, we prove that preserving all reduced decompositions is equivalent to the multiplier having a positive inverse. For a general DPO, the positive-inverse condition is equivalent both to flatness of the corresponding coefficient module and to the induced spatial map being a self-homeomorphism. Every strict extension of positive cones on the same field is nonflat, whereas preservation of a specified reduced decomposition is characterized by the vanishing of the corresponding first $\Tor$ obstruction.

To describe higher obstructions, we study topological realizations by the relevant complexes. If a lattice contains finite families of pairwise reduced nonzero positive elements of arbitrarily large size, its positive-sum complexes realize every finite simplicial complex; hence its Koszul homology may contain nonzero torsion. Under the same hypothesis, multiplicative truncation regions can have the homology of spheres of arbitrary dimension, producing nonzero multiplicative homology in arbitrarily high degrees.

For the concrete lattice of \cite{BY}, the coefficient modules of positive rational multipliers are flat, while those of all other nonzero positive multipliers have infinite flat dimension. We also obtain a Hahn lattice inequivalent to this example under field automorphisms and distinguish the two by their flat multipliers.

The paper has two parts. The first develops the homological language, order-carrying spaces, and comparison theorems. The second treats structural consequences, higher obstructions, and applications to concrete lattices. The calculations for the concrete lattice order on the complex field rely on the construction in \cite{BY}.

\part{Homology theory and order-carrying spaces}
\section{Homological obstructions to reduced decompositions}
\subsection{DPOs, lattice orders, and reduced decompositions}
Throughout, $L$ is a field of characteristic zero. We take $L=\C$ in applications to the complex field. Unless stated otherwise, positive cones need not contain $1$, be rationally saturated, or be finitely generated. All chain complexes are homologically indexed, and unspecified coefficients are $\Z$.

\begin{defn}[DPOs and reduced decompositions]\label{def:dpo}
A subset $P\subseteq L$ is a DPO positive cone if
\[
0\in P,\quad P+P\subseteq P,\quad PP\subseteq P,\quad
P\cap(-P)=\{0\},\quad P-P=L.
\]
Write $x\le_P y$ if and only if $y-x\in P$, and $p>_P0$ for $p\in P\setminus\{0\}$.
If this order is a lattice order, we call $P$ a lattice and write $x^+=0\vee x$ and $x^-=0\vee(-x)$. We call $P$ rationally saturated if $nx\in P$ implies $x\in P$ for every positive integer $n$.
A representation $x=p-q$, with $p,q\in P$, is reduced if every representation $x=a-b$, with $a,b\in P$, admits a unique expression
\[
a=p+c,\qquad b=q+c,\qquad c\in P.
\]
Positive elements $p,q$ are called reduced if $p-q$ is a reduced decomposition. Order comparisons without subscripts refer to $\le_P$; rational coefficients, real parameters, and heights are compared in the usual order.
\end{defn}

\begin{prop}[Reduced decompositions and least common upper bounds]\label{prop:reduced-min}
Let $V_x=P\cap(x+P)$. A representation $x=p-q$ is reduced if and only if $p$ is the least element of $V_x$. Thus every element has at most one reduced decomposition, and
\[
P\text{ is a lattice}\quad\Longleftrightarrow\quad
\text{every }x\in L\text{ has a reduced decomposition}.
\]
\end{prop}
\begin{proof}
See \byref{6}{Thm~2.9} and \byref{7}{Cor~2.10}. Their proofs use only additive closure and pointedness, not rational saturation.
\end{proof}

\begin{prop}[Basic lattice identities]\label{prop:lattice-basics}
Let $P$ be a lattice. Then
\[
x+y=(x\vee y)+(x\wedge y),\qquad
x=x^+-x^-,\quad x^+\wedge x^-=0,
\]
where $x^-=(-x)^+$. For $p,q\in P$, the representation $p-q$ is reduced if and only if $p\wedge q=0$. Moreover, $P$ is rationally saturated, and multiplication by nonnegative rational numbers preserves lattice operations.
\end{prop}
\begin{proof}
For the lattice identities and reducedness criterion, see \byref{6}{Thm~2.9}; for rational saturation, see \byref{5}{Lem~2.6}.
For $r\in\Q_{>0}$, we have $rP=P$, so $x\mapsto rx$ is an order automorphism and preserves $\vee,\wedge$. The conclusion also holds for $r=0$.
\end{proof}

\begin{cor}[Coordinate order on finite reduced families]\label{cor:disjoint-coordinates}
Let $P$ be a lattice. If $e_1,\ldots,e_m>_P0$ are pairwise reduced, then, for $r_i\in\Q$,
\[
\sum_i r_i e_i\in P\iff r_i\ge0\quad\text{for all }i.
\]
In particular, these elements are linearly independent, and $\bigvee_i e_i=\sum_i e_i$.
\end{cor}
\begin{proof}
If $u,v\ge0$ and $0\le z\le u+v$, then
\[
z=z_1+z_2,\qquad z_1=z\wedge u,\quad z_2=(z-u)^+,
\qquad 0\le z_1\le u,\quad0\le z_2\le v.
\]
Also $z_1,z_2\le z$. Hence $w\wedge u=w\wedge v=0$ implies $w\wedge(u+v)=0$.
For $r,s\in\Q_{\ge0}$, choose $M\ge r,s$. By Prop~\ref{prop:lattice-basics},
\[
0\le re_i\wedge se_j\le Me_i\wedge Me_j=M(e_i\wedge e_j)=0\quad(i\ne j).
\]
Thus $u=\sum_{r_i>0}r_ie_i$ and $v=\sum_{r_i<0}(-r_i)e_i$ satisfy $u\wedge v=0$.
If $u-v\in P$, then $v\le u$, whence $v=u\wedge v=0$. Pointedness excludes every negative coefficient. The converse follows from additive closure.
Applying this conclusion to $\sum r_ie_i=0$ and its negative gives linear independence. The last identity follows by induction from $u\vee v=u+v$ for reduced positive elements.
\end{proof}

\subsection{Linearization of positive translations}
\begin{defn}[Monoid ring and grading]\label{def:ring}
Set
\[
R=\Z[P]=\bigoplus_{p\in P}\Z X^p,
\qquad X^pX^q=X^{p+q},\qquad \deg X^p=p.
\]
Here $X^p$ is a formal symbol, and the identity of $R$ is $X^0$. Let $\Gr_LR$ denote the category of $L$-graded $R$-modules and degree-preserving maps. For a graded module $M$, use the convention $(M\sh v)_u=M_{u-v}$; in particular, $R\sh v$ is free of rank one with generator in degree $v$.
\end{defn}

\begin{prop}[Universal property of positive translations]\label{prop:translation}
For a fixed abelian group $M$, an $R$-module structure is equivalent to a monoid homomorphism $(P,+)\to(\operatorname{End}_{\Z}(M),\circ)$, $p\mapsto T_p$, that is,
\[
T_0=\id,\qquad T_pT_q=T_{p+q}.
\]
Moreover, multiplication by each $X^p$ on $R$ is injective, and monomial localization gives
\[
R[(X^p)^{-1}:p\in P]=\Z[L].
\]
Both rings are integral domains.
\end{prop}
\begin{proof}
Define the action by $X^pm=T_p(m)$; the relation $T_pT_q=T_{p+q}$ gives associativity. Conversely, take the actions of the monomials. These constructions are inverse.
On the free abelian basis $\{X^x:x\in L\}$, we have $X^pX^x=X^{p+x}$, so monomial multiplication is injective.
Since $P-P=L$,
\[
X^{p-q}=X^p(X^q)^{-1},\qquad R[(X^p)^{-1}:p\in P]=\Z[L].
\]
Choose a well-ordered $\Q$-basis of $L$ and define an additive total order $\preceq$ by the sign of the last nonzero coordinate of a nonzero vector.
If $f,g\in\Z[L]\setminus\{0\}$ have maximal support points $u,v$, the coefficient of $fg$ at $u+v$ is $f_ug_v\ne0$; every other pair of exponents has sum strictly less than $u+v$. Thus $\Z[L]$ and its subring $R$ have no zero divisors.
\end{proof}

\begin{cor}\label{cor:principal}
Each $X^pR$ is free, and $X^pR=X^qR$ if and only if $p=q$.
\end{cor}
\begin{proof}
Multiplication gives a graded isomorphism $R\sh p\to X^pR$. If $X^pR=X^qR$, comparison of monomial bases gives $p-q,q-p\in P$, hence $p=q$. The converse is immediate.
\end{proof}

\subsection{The binary complex and the reducedness criterion}
\begin{defn}[Binary obstruction to reducedness]\label{def:binary}
For $p,q\in P$, let $K(p,q)$ be the graded complex
\[
0\longrightarrow R\sh{p+q}
\xrightarrow{\ d_2\ }R\sh p\oplus R\sh q
\xrightarrow{\ d_1\ }R\longrightarrow0,
\]
where $d_2(r)=(-X^qr,X^pr)$ and $d_1(f,g)=X^pf+X^qg$. Write $\HH_P(p,q)=H_1K(p,q)$.
\end{defn}

\begin{thm}[Homological criterion for reducedness]\label{thm:binary}
We have $d_1d_2=0$ and natural isomorphisms
\[
\HH_P(p,q)\cong
\frac{X^pR\cap X^qR}{X^{p+q}R}
\cong\Tor_1^R(R/X^pR,R/X^qR).
\]
In particular, $\HH_P(p,q)=0$ if and only if $p-q$ is a reduced decomposition.
\end{thm}
\begin{proof}
We have $d_1d_2(r)=-X^{p+q}r+X^{q+p}r=0$. Set $I=X^pR\cap X^qR$. The maps
\[
\ker d_1\xrightarrow{\sim} I,\quad(f,g)\mapsto X^pf,
\qquad h\mapsto(h/X^p,-h/X^q),
\]
are inverse and send $\operatorname{im}d_2$ to $X^{p+q}R$.
Tensoring the free resolution $0\to R\sh p\xrightarrow{X^p}R\to R/X^pR\to0$ with $R/X^qR$, its first kernel is identified with $I/X^{p+q}R$ by $[f]\mapsto[X^pf]$.
By the monomial bases, $I/X^{p+q}R=0$ is equivalent to
\[
(p+P)\cap(q+P)=p+q+P.
\]
The bijection $(a,b)\mapsto a+q=b+p$ identifies positive representations $p-q=a-b$ with the left-hand side. Its value is $p+q+c$ if and only if $(a,b)=(p+c,q+c)$, which is precisely reducedness.
\end{proof}

\begin{cor}[Representations with vanishing obstruction]\label{cor:zero-pairs}
Let $\mathcal Z_P=\{(p,q)\in P^2:\HH_P(p,q)=0\}$. The map
\[
\mathcal Z_P\longrightarrow L,\qquad(p,q)\longmapsto p-q
\]
is injective, and it is surjective if and only if $P$ is a lattice.
\end{cor}
\begin{proof}
Thm~\ref{thm:binary} identifies $\mathcal Z_P$ with all reduced representations. Uniqueness in Prop~\ref{prop:reduced-min} gives injectivity, and its existence criterion gives the assertion about surjectivity.
\end{proof}

\subsection{Homological recovery of common positive parts}
\begin{thm}[Recovery of the common positive part]\label{thm:meet-recovery}
Let $P$ be a lattice, $c=p\wedge q$, and $b=p\vee q$. Then, as graded modules,
\[
\HH_P(p,q)\cong X^bR/X^{b+c}R
\cong(R/X^cR)\sh b,
\qquad \Ann_R\HH_P(p,q)=X^cR.
\]
Thus the annihilator of this homology module uniquely determines $p\wedge q$.
\end{thm}
\begin{proof}
The lattice order gives $(p+P)\cap(q+P)=b+P$ and $p+q=b+c$, so Thm~\ref{thm:binary} yields the asserted isomorphisms.
For any ideal $I$, we have $r(R/I)=0\iff r\in I$, hence $\Ann_R(R/X^cR)=X^cR$. Uniqueness of $c$ follows from Cor~\ref{cor:principal}.
\end{proof}

\begin{cor}\label{cor:fixed-pair}
Let $P$ be a lattice, $p,q\in P$, and $c=p\wedge q$. Then
\[
p-q=(p-c)-(q-c)
\]
is a reduced decomposition, and $\HH_P(p,q)\ne0$ if and only if $c\ne0$.
\end{cor}
\begin{proof}
Since $(p-c)\wedge(q-c)=(p\wedge q)-c=0$, the pair $(p-c,q-c)$ is reduced. Thm~\ref{thm:meet-recovery} gives $\HH_P(p,q)\ne0\iff c\ne0$.
\end{proof}

\section{Koszul homology and positive multiplication}
\subsection{Multivariable Koszul homology}
\begin{defn}[Multivariable complex]\label{def:koszul}
For $m\ge1$ and $\mathbf p=(p_1,\ldots,p_m)\in P^m$, write $[m]=\{1,\ldots,m\}$ and set
\[
K_j(\mathbf p)=\bigoplus_{I\subseteq[m],\,|I|=j}R\sh{p_I},\qquad
p_I=\sum_{i\in I}p_i.
\]
Let $e_I$ denote the corresponding generator. For $I=\{i_1<\cdots<i_j\}$, define
\[
d e_I=\sum_{r=1}^j(-1)^{r-1}X^{p_{i_r}}e_{I\setminus\{i_r\}}.
\]
\end{defn}

\begin{prop}\label{prop:koszul-d}
These maps preserve degree and satisfy $d^2=0$. Moreover, $H_0K(\mathbf p)=R/(X^{p_1},\ldots,X^{p_m})$, and $H_mK(\mathbf p)=0$ for $m\ge1$.
\end{prop}
\begin{proof}
For $I=\{i_1<\cdots<i_j\}$, each summand has degree $p_{i_r}+p_{I\setminus\{i_r\}}=p_I$. For $r<s$, the coefficient obtained by deleting $i_r,i_s$ in $d^2e_I$ is
\[
\bigl((-1)^{r+s-3}+(-1)^{r+s-2}\bigr)X^{p_{i_r}+p_{i_s}}=0.
\]
The identity $\operatorname{im}d_1=(X^{p_1},\ldots,X^{p_m})$ gives $H_0$. If $d_m(fe_{[m]})=0$, one component is $\pm X^{p_1}f=0$, so $f=0$ and hence $H_m=0$.
\end{proof}

\begin{thm}[Top relations and the common positive part]\label{thm:top-koszul}
Let $P$ be a lattice, $m\ge2$, $c=\bigwedge_{i=1}^m p_i$, and $s=\sum_i p_i$. Then
\[
H_{m-1}K(\mathbf p)\cong(R/X^cR)\sh{s-c}.
\]
\end{thm}
\begin{proof}
Set $S=\Z[L]$. On $K(\mathbf p)\otimes_RS$,
\[
h(\xi)=X^{-p_1}e_1\wedge\xi,\qquad dh+hd=\id.
\]
Thus every cycle in degree $m-1$ is $f\,de_{[m]}$, with $f\in S$ unique because $d_m$ is injective.
This cycle belongs to $K_{m-1}(\mathbf p)$ if and only if $fX^{p_i}\in R$ for every $i$, that is,
\[
f\in\bigcap_iX^{-p_i}R=X^{-c}R;
\qquad w+p_i\in P\ (\forall i)\iff w\ge_P\bigvee_i(-p_i)=-c.
\]
Boundaries correspond to $f\in R$, and $f\mapsto f\,de_{[m]}$ increases degree by $s$. Hence
\[
H_{m-1}K(\mathbf p)\cong(X^{-c}R/R)\sh s
\cong(R/X^cR)\sh{s-c}.
\]
\end{proof}

\begin{cor}\label{cor:top-mixing}
Let $P$ be a lattice and $m\ge2$. If $p_1,\ldots,p_m$ are pairwise reduced and $a\in P\setminus\{0\}$, then
\[
H_{m-1}K(ap_1,\ldots,ap_m)\ne0
\iff \bigwedge_i ap_i>_P0.
\]
\end{cor}
\begin{proof}
Set $c=\bigwedge_iap_i\in P$. By Thm~\ref{thm:top-koszul}, this homology is nonzero if and only if $X^cR\ne R$, which is equivalent to $c\ne0$ by Cor~\ref{cor:principal}.
\end{proof}

\subsection{Finite geometric models for graded homology}
\begin{defn}[Positive-sum compatibility complex]\label{def:delta}
For $u\in P$, define the finite simplicial complex
\[
\Delta_u(\mathbf p)=\{I\subseteq[m]:p_I\le_Pu\}.
\]
It contains the empty face; some indices may not be vertices. Reduced chains use the augmented convention: the complex consisting only of the empty face has homology $\Z$ in degree $-1$.
\end{defn}

\begin{thm}[Graded chain comparison]\label{thm:delta}
There are natural identifications of chain complexes
\[
K_j(\mathbf p)_u\cong\widetilde C_{j-1}(\Delta_u(\mathbf p);\Z),
\]
and consequently, for $j\ge0$,
\[
H_jK(\mathbf p)_u\cong\widetilde H_{j-1}(\Delta_u(\mathbf p);\Z).
\]
For $r\in P$, the module action of $X^r$ corresponds to the map induced by the inclusion $\Delta_u\subseteq\Delta_{u+r}$.
\end{thm}
\begin{proof}
Define $\theta_{j,u}(X^{u-p_I}e_I)=[i_1,\ldots,i_j]$ for $I\in\Delta_u$; for $j=0$, send $X^u$ to $1$ in the augmentation term. These maps are degreewise bijections of bases, and
\[
\theta_{j-1,u}d(X^{u-p_I}e_I)
=\sum_{r=1}^j(-1)^{r-1}[I\setminus\{i_r\}]
=\partial\theta_{j,u}(X^{u-p_I}e_I).
\]
Moreover, $\theta_{j,u+r}(X^rX^{u-p_I}e_I)=[I]$, so they commute with the maps induced by face inclusions. Passing to homology proves the claim.
\end{proof}

\begin{cor}\label{cor:relative}
For $u,r\in P$, let $A=\Delta_u(\mathbf p)$ and $B=\Delta_{u+r}(\mathbf p)$. There is a natural long exact sequence
\[
\cdots\to H_j(B,A;\Z)\xrightarrow{\partial}H_jK(\mathbf p)_u
\xrightarrow{X^r}H_jK(\mathbf p)_{u+r}
\to H_{j-1}(B,A;\Z)\to\cdots .
\]
where $j\ge0$ and relative homology in negative degrees is zero. In particular, $\ker(X^r)=\operatorname{im}\partial$.
\end{cor}
\begin{proof}
Apply the long exact homology sequence to
\[
0\to\widetilde C_\bullet(A;\Z)\to\widetilde C_\bullet(B;\Z)
\to C_\bullet(B,A;\Z)\to0
\]
and use the degree shift and naturality in Thm~\ref{thm:delta}.
\end{proof}

\subsection{Region ideals and their resolutions}
\begin{prop}[Dependence on the generating family]\label{prop:redundancy}
For $p>_P0$, the ideals $(X^p)=(X^p,X^{2p})$ coincide, but
\[
H_1K(p)=0,\qquad H_1K(p,2p)\cong(R/X^pR)\sh{2p}.
\]
Thus Koszul homology records relations among the specified generators.
\end{prop}
\begin{proof}
The identity $X^{2p}=X^pX^p$ gives equality of the ideals. By Prop~\ref{prop:koszul-d}, $H_1K(p)=0$; by Thm~\ref{thm:binary},
\[
H_1K(p,2p)=X^{2p}R/X^{3p}R\cong(R/X^pR)\sh{2p}\ne0,
\]
The last step uses $p>0$ and Cor~\ref{cor:principal}.
\end{proof}

\begin{defn}[Region ideal]\label{def:monomial}
For $V\subseteq P$ with $V+P\subseteq V$, set $J_V=\bigoplus_{v\in V}\Z X^v$.
If $V=\bigcup_i(p_i+P)$, then $J_V=\sum_i X^{p_i}R$. In what follows, fix $m\ge1$ and $p_1,\ldots,p_m\in P$, and write $J_i=X^{p_i}R$ and $J=\sum_{i=1}^mJ_i$.
\end{defn}

\begin{prop}[Exact intersection complex]\label{prop:intersection-resolution}
For any DPO positive cone, signed inclusions give an exact complex
\[
0\longrightarrow\bigcap_iJ_i\longrightarrow\cdots\longrightarrow
\bigoplus_{|I|=2}\bigcap_{i\in I}J_i\longrightarrow
\bigoplus_iJ_i\longrightarrow R\longrightarrow R/J\longrightarrow0.
\]
\end{prop}
\begin{proof}
Fix $u\in L$ and let $A_u=\{i:u-p_i\in P\}$. Each intersection module has basis $X^u$ in degree $u$. If $A_u\ne\varnothing$, this degree component is the augmented chain complex of the simplex on $A_u$; coning to any vertex gives a chain contraction.
If $A_u=\varnothing$, all intersection terms vanish, and the remaining map $R_u\to(R/J)_u$ is the identity. Thus the complex is exact in every degree.
\end{proof}

\begin{thm}[Finite free model for a lattice]\label{thm:taylor}
Let $P$ be a lattice, $b_I=\bigvee_{i\in I}p_i$, and $b_\varnothing=0$. Then $R/J$ has a free resolution
\[
T_j=\bigoplus_{|I|=j}R\sh{b_I},\qquad
de_I=\sum_{r=1}^{j}(-1)^{r-1}X^{b_I-b_{I\setminus\{i_r\}}}e_{I\setminus\{i_r\}}.
\]
It admits a differential graded algebra structure. For $I\cap I'=\varnothing$, let $\eps(I,I')$ be the sign of the permutation that orders the concatenation $(I,I')$ increasingly, and set
\[
e_Ie_{I'}=\eps(I,I')X^{b_I+b_{I'}-b_{I\cup I'}}e_{I\cup I'}
=\eps(I,I')X^{b_I\wedge b_{I'}}e_{I\cup I'};
\]
For $I\cap I'\ne\varnothing$, set $e_Ie_{I'}=0$.
\end{thm}
\begin{proof}
For $I\ne\varnothing$, we have $\bigcap_{i\in I}X^{p_i}R=X^{b_I}R$. Under $e_I\mapsto X^{b_I}$, the differential becomes the intersection inclusions of Prop~\ref{prop:intersection-resolution}, so the complex resolves $R/J$.
Over $S=\Z[L]$, denote the Koszul basis by $\epsilon_I$ and define a degreewise injection
\[
\theta:T_\bullet\hookrightarrow K(\mathbf p)\otimes_RS,
\qquad\theta(e_I)=X^{b_I-p_I}\epsilon_I.
\]
Direct substitution gives $\theta d=d\theta$ and $\theta(e_Ie_{I'})=\theta(e_I)\wedge\theta(e_{I'})$.
Since $b_I+b_{I'}-b_{I\cup I'}=b_I\wedge b_{I'}\in P$, the product lies in $T$. By injectivity, associativity, graded commutativity, and
\[
d(xy)=(dx)y+(-1)^{\deg_{\rm hom}x}x(dy)
\]
are inherited from the Koszul exterior algebra.
\end{proof}

\begin{cor}\label{cor:independent}
Let $P$ be a lattice, and let two finite families of positive elements generate the same monomial ideal $J$. The graded free resolutions of Thm~\ref{thm:taylor} are chain homotopy equivalent. Hence, for every $R$-module $M$ and $n\ge0$, they compute naturally isomorphic groups $\Tor_n^R(R/J,M)$ and $\Ext_R^n(R/J,M)$.
\end{cor}
\begin{proof}
If $T,T'$ resolve the same $R/J$, the comparison theorem for projective resolutions gives degree-preserving chain maps $f:T\to T'$ and $g:T'\to T$ lifting the identity, with
\[
gf-\id=dh+hd,\qquad fg-\id=dh'+h'd.
\]
Tensoring or applying $\Hom$ preserves these chain homotopy identities and hence induces inverse isomorphisms on homology. Different lifts are chain homotopic, so the isomorphisms are natural and independent of the lifts.
\end{proof}

\subsection{Positive multiplication and change of coefficients}
For $a\in P\setminus\{0\}$, positive multiplication $x\mapsto ax$ preserves order but need not preserve reduced decompositions. We describe this obstruction by homology after a change of coefficients.

\begin{defn}[Multiplicative change of coefficients]\label{def:twist}
For $a\in P\setminus\{0\}$, let $\phi_a(X^p)=X^{ap}$ and write
\[
R_a={}_{\phi_a}R_R.
\]
Its left action is $X^p\cdot f=X^{ap}f$, and its right action is ordinary multiplication. Flatness of $R_a$ refers to its left $R$-module structure. Viewing $M$ as a right $R$-module, $M\otimes_RR_a$ uses this left action; the right $R$-actions on the tensor product and its $\Tor$ groups are induced by the right action on $R_a$. For homogeneous elements, set
\[
\deg(m\otimes f)=a\deg m+\deg f;
\]
Since $a(v+p)+w=av+(ap+w)$, this grading respects the tensor relations. In particular, a free generator originally in degree $v$ acquires degree $av$.
\end{defn}

\begin{thm}[Binary multiplicative defect]\label{thm:defect}
Let $P$ be a lattice, $p,q\in P$, and $a\in P\setminus\{0\}$. Set
\[
\delta_a(p,q)=a(p\vee q)-(ap\vee aq)
=ap\wedge aq-a(p\wedge q)\in P.
\]
Then
\[
\Tor_1^R(R/(X^p,X^q),R_a)
\cong\frac{X^{ap\vee aq}R}{X^{a(p\vee q)}R}
\cong(R/X^{\delta_a(p,q)}R)\sh{ap\vee aq}.
\]
\end{thm}
\begin{proof}
Set $b=p\vee q$, $c=p\wedge q$, $b'=ap\vee aq$, and $c'=ap\wedge aq$. Since $ab\ge ap,aq$ and $b+c=p+q$,
\[
ab-b'\in P,\qquad ab-b'=a(p+q-c)-(ap+aq-c')=c'-ac.
\]
After tensoring the binary Taylor resolution with $R_a$, the differentials are
\[
d_1(f,g)=X^{ap}f+X^{aq}g,\qquad
d_2(h)=(-X^{a(b-p)}h,X^{a(b-q)}h).
\]
The isomorphism $\ker d_1\to X^{ap}R\cap X^{aq}R$, $(f,g)\mapsto X^{ap}f$, sends boundaries to $X^{ab}R$. Hence
\[
\Tor_1^R(R/J,R_a)=X^{b'}R/X^{ab}R
\cong(R/X^{ab-b'}R)\sh{b'}.
\]
\end{proof}

\begin{defn}[Multiplicative compatibility complex]\label{def:xi}
In the lattice case, for $u\in P$, define
\[
\Xi_{a,u}(\mathbf p)=\{I\subseteq[m]:a b_I\le_Pu\}.
\]
\end{defn}

\begin{thm}[Finite computation of multiplicative homology]\label{thm:xi}
For $n\ge1$,
\[
\Tor_n^R(R/J,R_a)_u
\cong\widetilde H_{n-1}(\Xi_{a,u}(\mathbf p);\Z).
\]
This isomorphism identifies the action of $X^r$ with the map induced by inclusion as $u$ changes to $u+r$.
\end{thm}
\begin{proof}
We have $T_j\otimes_RR_a=\bigoplus_{|I|=j}R\sh{ab_I}$. In degree $u$, define
\[
X^{u-ab_I}e_I\longmapsto[I],\qquad I\in\Xi_{a,u}(\mathbf p).
\]
Since $X^{a(b_I-b_{I\setminus\{i\}})}X^{u-ab_I}=X^{u-ab_{I\setminus\{i\}}}$, this basis bijection identifies the differential with the augmented boundary and hence gives the homology isomorphism. Multiplication by $X^r$ preserves $[I]$ and therefore corresponds to inclusion of complexes.
\end{proof}

\begin{cor}\label{cor:binary-high-vanish}
For a lattice, $\Tor_n^R(R/(X^p,X^q),R_a)=0$ for every $n\ge2$.
\end{cor}
\begin{proof}
The binary Taylor resolution vanishes in degrees $>2$. After changing coefficients, one component of $d_2$ is multiplication by the nonzero monomial $X^{a((p\vee q)-p)}$, so $d_2$ is injective. Hence $\Tor_2=0$, and all higher degrees vanish as well.
\end{proof}

\section{Order-carrying spaces and obstructions to the existence of reduced decompositions}
\subsection{Order-carrying spaces induced by positive cones}
\begin{defn}[Order-carrying topology]\label{def:carrier}
For a DPO positive cone $P$, define
\[
\tau_P=\{U\subseteq L:U+P\subseteq U\},\qquad X_P=(L,\tau_P),
\qquad U_x=x+P.
\]
The topology is always the upper Alexandrov topology; it generally differs from the Euclidean topology when $L=\C$.
\end{defn}

\begin{prop}[Basic properties]\label{prop:carrier-basics}
The space $X_P$ is a $T_0$ Alexandrov space, and $U_x$ is the minimal open neighborhood of $x$. Translations are self-homeomorphisms, addition is continuous, and each positive multiplication $m_a(x)=ax$ is continuous. For every $A\subseteq L$,
\[
\overline A=A-P,\qquad
\operatorname{Int}A=\{x:x+P\subseteq A\}.
\]
Any two nonempty open sets intersect, and
\[
P=U_0,\qquad
x\le_Py\iff y\in U_x,\qquad
P\text{ is a lattice}\iff\forall x,y\ \exists z:\ U_x\cap U_y=U_z.
\]
\end{prop}
\begin{proof}
The family $\tau_P$ is closed under arbitrary unions and intersections. If $x\in O\in\tau_P$, then $U_x\subseteq O$, so $U_x$ is minimal. If $U_x=U_y$, then $x-y,y-x\in P$; pointedness gives $x=y$, proving the $T_0$ property.
Let $t_c(x)=x+c$. Then
\[
t_c(U_x)=U_{x+c},\qquad U_x+U_y=U_{x+y},\qquad
m_a(U_x)\subseteq U_{ax}\quad(a\in P).
\]
These equalities and inclusions give, respectively, homeomorphisms under translation, continuity of addition, and continuity of positive multiplication. Also,
\[
x\in\overline A\iff U_x\cap A\ne\varnothing\iff x\in A-P,
\qquad x\in\operatorname{Int}A\iff U_x\subseteq A.
\]
Write $y-x=r-s$ with $r,s\in P$. Then $x+r=y+s\in U_x\cap U_y$, so any two nonempty open sets intersect.
Finally, $z$ is the least element of $U_x\cap U_y$ if and only if this set equals $U_z$, equivalently, $z=x\vee y$. Infima follow by negation.
\end{proof}

\subsection{Graded modules and sheaves}
\begin{defn}[Sheaf associated to a graded module]\label{def:sheaf}
For $M\in\Gr_LR$, the groups $M_u$ and the maps
\[
M_u\xrightarrow{X^{v-u}}M_v\qquad(u\le_Pv)
\]
define a sheaf $\widetilde M$. Explicitly, for an open set $U$, $\widetilde M(U)$ consists of compatible families $(s_u)_{u\in U}$ satisfying $X^{v-u}s_u=s_v$ whenever $u\le v$. Its stalk at $u$ is $M_u$.
\end{defn}

\begin{thm}[Equivalence of graded modules and sheaves]\label{thm:sheaf-equivalence}
This construction gives an exact equivalence of categories
\[
\Gr_LR\simeq\Sh(X_P;\Ab).
\]
For an upward-closed subset $V\subseteq P$, $M\in\Gr_LR$, and $n\ge0$, there are natural isomorphisms
\[
\Ext^n_{\Gr_LR}(J_V,M)
\cong H^n(V;\widetilde M|_V).
\]
The left-hand side uses only degree-preserving maps; the right-hand side is ordinary sheaf cohomology.
\end{thm}
\begin{proof}
For a sheaf $\mathcal F$, define
\[
G(\mathcal F)=\bigoplus_{u\in L}\mathcal F(U_u),\qquad
X^p s=\operatorname{res}_{U_u,U_{u+p}}(s).
\]
Transitivity of restriction maps gives the $R$-module action. For $M\in\Gr_LR$, the evaluation maps
\[
\widetilde M(U_u)\xrightarrow{\sim}M_u,
\qquad (s_v)\mapsto s_u,
\quad m\mapsto(X^{v-u}m)_{v\in U_u}
\]
are inverse isomorphisms. For an open set $O$, the sheaf gluing axiom identifies $\mathcal F(O)$ with
\[
\{(s_u)\in\prod_{u\in O}\mathcal F(U_u):
\operatorname{res}_{U_u,U_v}s_u=s_v\text{ whenever }u\le v\},
\]
so $G$ and $M\mapsto\widetilde M$ are naturally inverse. Since the stalk at $u$ equals $\mathcal F(U_u)$, exactness is tested degreewise or stalkwise on the respective sides.
A degree-preserving map $f:J_V\to M$ is determined by $s_u=f(X^u)\in M_u$, and
\[
X^{v-u}s_u=s_v\quad(u\le v),\qquad
\Hom_{\Gr_LR}(J_V,M)\cong\Gamma(V,\widetilde M|_V).
\]
For the open embedding $j:V\hookrightarrow X_P$, extension by zero $j_!$ is stalkwise exact and satisfies $j_!\dashv j^*$; hence $j^*$ preserves injectives. If $M\to I^\bullet$ is a graded injective resolution, then
\[
H^n\Hom_{\Gr_LR}(J_V,I^\bullet)
=H^n\Gamma(V,\widetilde I^\bullet|_V)
=H^n(V;\widetilde M|_V).
\]
\end{proof}

\begin{cor}\label{cor:region-cohomology}
For $x\in L$, set $V_x=P\cap(x+P)$ and $J_x=J_{V_x}$. For every $M\in\Gr_LR$ and $n\ge0$, there are natural isomorphisms
\[
\Ext^n_{\Gr_LR}(J_x,M)\cong H^n(V_x;\widetilde M|_{V_x}).
\]
\end{cor}
\begin{proof}
Since $V_x\subseteq P$ and $V_x+P\subseteq V_x$, apply Thm~\ref{thm:sheaf-equivalence} with $V=V_x$.
\end{proof}

\subsection{Canonical obstructions to the existence of reduced decompositions}
\begin{defn}[Sheaf of compatible choices and obstruction class]\label{def:obstruction}
On $V_x$, set
\[
(\cE_x)_u=\bigoplus_{\substack{t\in V_x\\t\le_Pu}}\Z e_t,
\]
with inclusions as transition maps. Summing coefficients gives a surjection $\cE_x\to\underline\Z$; denote its kernel by $\cN_x$. Define
\[
o_P(x)=\delta(1)\in H^1(V_x;\cN_x),
\]
where $\delta$ comes from the long exact cohomology sequence of $0\to\cN_x\to\cE_x\to\underline\Z\to0$.
\end{defn}

\begin{thm}[Canonical obstruction to the existence of a reduced decomposition]\label{thm:obstruction}
For every DPO positive cone $P$ and $x\in L$,
\[
o_P(x)=0\iff V_x\text{ has a least element}\iff x\text{ has a reduced decomposition}.
\]
\end{thm}
\begin{proof}
On each stalk, $e_u\mapsto1$, so the defining sequence is exact. Its long exact sequence of global sections gives
\[
o_P(x)=0\iff1\in\operatorname{im}\bigl(\Gamma(V_x,\cE_x)\to\Gamma(V_x,\underline\Z)\bigr).
\]
Set $F=\bigoplus_{t\in V_x}\Z e_t$. For a compatible section $(s_u)$, any $u,v\in V_x$ have a common upper bound $w\in V_x$. Since all transition maps are inclusions, $s_u=s_w=s_v$ in $F$. Thus $(s_u)$ corresponds to a single finitely supported vector $s=\sum_tn_te_t$, with
\[
\supp s\subseteq\bigcap_{u\in V_x}\{t\in V_x:t\le u\}.
\]
If $\sum_tn_t=1$, its support is nonempty, and any support point is the least element of $V_x$. Conversely, a least element $t$ gives the section $s_u=e_t$ mapping to $1$. Apply Prop~\ref{prop:reduced-min} for the last equivalence.
\end{proof}

\begin{cor}[Cohomological criterion for lattice ordering]\label{cor:lattice-obstruction}
The cone $P$ is a lattice if and only if $o_P(x)=0$ for all $x\in L$.
\end{cor}
\begin{proof}
By Thm~\ref{thm:obstruction}, all $o_P(x)$ vanish if and only if every $x$ has a reduced decomposition. By Prop~\ref{prop:reduced-min}, this is equivalent to $P$ being a lattice.
\end{proof}

\subsection{Reducedness obstructions, least elements, and compactness}
Here compactness means that every open cover has a finite subcover; no Hausdorff condition is imposed.

\begin{prop}[Finite-cover criterion for compactness]\label{prop:compact-cover}
A subset $K\subseteq X_P$ is compact if and only if there is a finite set $A\subseteq K$ such that $K\subseteq A+P$. In particular, every $U_x$ is compact, and an upward-closed set $U$ is compact if and only if it is a finite union of sets $U_a$ with $a\in U$.
\end{prop}
\begin{proof}
If $K$ is compact, take a finite subcover of $\{U_x:x\in K\}$ to obtain the required $A\subseteq K$.
Conversely, given an open cover $K\subseteq\bigcup_iO_i$, choose $i(a)$ with $a\in O_{i(a)}$ for each $a\in A$. Minimality of the neighborhoods gives
\[
K\subseteq\bigcup_{a\in A}U_a\subseteq\bigcup_{a\in A}O_{i(a)},
\]
so $K$ is compact. If $K=U$ is upward closed, then $A\subseteq U$ implies $A+P\subseteq U$, hence $U=A+P$.
\end{proof}

\begin{thm}[Compactness and least common upper bounds]\label{thm:compact-min}
Suppose $P$ is rationally saturated. For every $x,y\in L$,
\[
U_x\cap U_y\text{ is compact}\iff x\vee y\text{ exists}.
\]
In particular, $o_P(x)=0$ if and only if $V_x$ is compact.
\end{thm}
\begin{proof}
First choose an additive total order extending $P$. Explicitly, set $\operatorname{Sat}(P)=\{x\in L:nx\in P\text{ for some integer }n\ge1\}$, a pointed rational cone by \byref{5}{Lem~2.5}. By Zorn's lemma, choose a maximal pointed rational cone $C$ containing it. If $v,-v\notin C$, then $C+\Q_{\ge0}v$ is still pointed: if
\[
c+rv=-(c'+r'v),\qquad c,c'\in C,\quad r,r'\ge0,
\]
then $r+r'>0$ would give $-v=(c+c')/(r+r')\in C$, whereas $r+r'=0$ gives $c=c'=0$. This contradicts maximality, so $C\cup(-C)=L$ and defines the required additive total order $\preceq$.
If $z=x\vee y$ exists, then $U_x\cap U_y=U_z$ is compact. Conversely, by Prop~\ref{prop:compact-cover}, write
\[
U_x\cap U_y=\bigcup_{j=1}^m U_{z_j},
\]
and remove redundant terms so that the $z_j$ are pairwise incomparable. If $m\ge2$, let the first two terms in the $\preceq$-ordering be $z_1\prec z_2$ and set $w=(z_1+z_2)/2$. Then
\[
w-x=\tfrac12(z_1-x+z_2-x)\in P,\qquad w-y\in P.
\]
so $z_j\le_Pw$ for some $j$. Since $z_1\prec w\prec z_2$ and $P\subseteq C$, necessarily $j=1$, whence $z_2-z_1=2(w-z_1)\in P$, a contradiction. Thus $m=1$ and $z_1=x\vee y$. Finally, take $(x,y)=(0,x)$ and apply Thm~\ref{thm:obstruction}.
\end{proof}

\begin{defn}[Coherence]\label{def:coherence}
The space $X_P$ is coherent if the intersection of any two compact saturated sets is compact. Here saturated sets are intersections of open sets; in an Alexandrov space, they are precisely the open sets.
\end{defn}

\begin{cor}[Compactness criterion for lattice ordering]\label{cor:coherence}
A rationally saturated DPO positive cone $P$ is a lattice if and only if its order-carrying space is coherent, equivalently, if every $U_x\cap U_y$ is compact.
\end{cor}
\begin{proof}
If every $U_x\cap U_y$ is compact, then, for compact open sets $A=\bigcup_{i=1}^mU_{x_i}$ and $B=\bigcup_{j=1}^nU_{y_j}$,
\[
A\cap B=\bigcup_{i,j}(U_{x_i}\cap U_{y_j}),
\]
is a finite union of compact sets and hence compact. Conversely, take $A=U_x,B=U_y$. Now apply Thm~\ref{thm:compact-min} and Prop~\ref{prop:carrier-basics}.
\end{proof}

\subsection{The necessity of coefficient sheaves}
\begin{defn}[Order complex]\label{def:order-complex}
For a poset $A$, let $\ord(A)$ be its order complex: the nonempty simplices are the finite nonempty strict chains, and the empty face is included by convention. Its geometric realization carries the weak topology.
\end{defn}

\begin{prop}[Vanishing with constant coefficients on directed systems]\label{prop:directed-acyclic}
If $A$ is nonempty and upward directed, then $\widetilde H_n(\ord(A);\Z)=0$ for every $n\ge0$. In particular, this holds for the order complexes of $L$, every nonempty open set, and every $V_x$.
\end{prop}
\begin{proof}
An augmented chain cycle $z$ involves only finitely many vertices. Choose a common upper bound $b$ of these vertices and, on the subcomplex they generate together with $b$, set
\[
h_k[v_0,\ldots,v_k]=
\begin{cases}
(-1)^{k+1}[v_0,\ldots,v_k,b],&v_k<b,\\
0,&v_k=b,
\end{cases}
\qquad h_{-1}(1)=[b].
\]
The face-deletion formula gives $\partial h+h\partial=\id$, so $z=\partial h(z)$.
The set $L$ is upward directed. Nonempty open sets and $V_x$ are upward closed, so common upper bounds of their finite subsets remain in them.
\end{proof}

\begin{cor}\label{cor:coeff-necessary}
Let $P$ be a DPO positive cone. For every $x\in L$ and $n\ge0$,
\[
\widetilde H_n(\ord(V_x);\Z)=0,
\qquad
o_P(x)=0\iff V_x\text{ has a least element}.
\]
\end{cor}
\begin{proof}
Apply Prop~\ref{prop:directed-acyclic} to the upward-directed set $V_x$ to obtain $\widetilde H_n(\ord(V_x);\Z)=0$ for $n\ge0$. Meanwhile, Thm~\ref{thm:obstruction} gives $o_P(x)=0\iff V_x$ has a least element.
\end{proof}

\section{Common-upper-bound modules and criteria for lattice ordering}
\subsection{Common-upper-bound regions and their modules}
\begin{defn}[Common-upper-bound module]\label{def:upper-module}
For $p,q\in P$, write
\[
U(p,q)=(p+P)\cap(q+P),\qquad J(p,q)=J_{U(p,q)}.
\]
For any nonempty upward-closed subset $V\subseteq P$, define the canonical presentation
\[
E_V=\bigoplus_{t\in V}R\sh t\xrightarrow{\eps_V}J_V,
\qquad e_t\longmapsto X^t,
\]
Write $N_V=\ker\eps_V$, and let $\omega_V\in\Ext^1_{\Gr_LR}(J_V,N_V)$ be the extension class of this short exact sequence.
\end{defn}

\begin{prop}\label{prop:upper-system}
The set $U(p,q)$ is nonempty and upward closed, and
\[
U(p,q)=q+V_{p-q}.
\]
Every $J_V$ is the sum of its rank-one free monomial submodules $X^tR$, $t\in V$. It is a principal monomial ideal if and only if $V$ has a least element.
\end{prop}
\begin{proof}
We have $p+q\in U(p,q)$, and
\[
v\in U(p,q)\iff v-q\in P\cap(p-q+P)=V_{p-q}.
\]
The inclusion $V+P\subseteq V$ gives $\sum_{t\in V}X^tR\subseteq J_V$. The reverse inclusion holds because each basis element satisfies $X^v\in X^vR$.
Comparing monomial bases gives $J_V=X^tR\iff V=t+P$, the latter being equivalent to $t=\min V$.
\end{proof}

\subsection{Downward directedness and least elements}
\begin{thm}[Projectivity and least elements]\label{thm:projective}
For a nonempty upward-closed subset $V\subseteq P$, the following are equivalent:
\begin{enumerate}
\item $J_V$ is projective in $\Gr_LR$;
\item $\omega_V=0$;
\item $V$ has a least element.
\end{enumerate}
For $V=V_x$, Thm~\ref{thm:sheaf-equivalence} identifies $\omega_V$ with $o_P(x)$.
\end{thm}
\begin{proof}
Since $E_V$ is graded free, $J_V$ is projective if and only if the canonical short exact sequence splits, equivalently, $\omega_V=0$. If $s:J_V\to E_V$ is a degree-preserving splitting, write
\[
s(X^u)=\sum_{t\le u}n_t(u)X^{u-t}e_t,\qquad\sum_tn_t(u)=1.
\]
For $u\le w$, the identity $s(X^w)=X^{w-u}s(X^u)$ gives $n_t(w)=n_t(u)$. Any $u,v\in V$ have a common upper bound, so all coefficient vectors equal a single finite vector $(n_t)$. Any $t$ in its nonempty support satisfies $t\le u$ for all $u\in V$, hence $t=\min V$.
Conversely, if $t=\min V$, then $J_V=X^tR\cong R\sh t$ is projective.
For $V=V_x$, we have $\widetilde J_V|_V=\underline\Z$, and the canonical short exact sequence restricts to $0\to\cN_x\to\cE_x\to\underline\Z\to0$. Under the isomorphism of Thm~\ref{thm:sheaf-equivalence}, the connecting map sends the section $1$ corresponding to $\id_{J_V}$ to $\omega_V$; its image is therefore $o_P(x)$.
\end{proof}

\begin{thm}[Flatness and downward directedness]\label{thm:flat-ideal}
For a nonempty upward-closed subset $V\subseteq P$, $J_V$ is flat as an ordinary $R$-module if and only if $V$ is downward directed, that is,
\[
\forall u,v\in V\quad\exists w\in V:\quad w\le_Pu,v.
\]
\end{thm}
\begin{proof}
If $V$ is downward directed, then $\{X^wR:w\in V\}$ is directed under inclusion, and
\[
J_V=\varinjlim_{w\in V}X^wR.
\]
Each term is free, and directed colimits preserve flatness; hence $J_V$ is flat.
Conversely, take $u,v\in V$ and a common lower bound $c\in L$. Set $r=X^{u-c}$ and $s=X^{v-c}$, so $sX^u-rX^v=0$. The equational criterion for flatness \cite[Lem~10.39.11]{Stacks} gives finitely many $z_\ell\in J_V$ and $A_\ell,B_\ell\in R$ such that
\[
X^u=\sum_\ell A_\ell z_\ell,\qquad
X^v=\sum_\ell B_\ell z_\ell,\qquad sA_\ell=rB_\ell.
\]
Writing $[X^x]f$ for the coefficient of $X^x$ in $f$, we have
\[
1=\sum_{\ell,w}[X^{u-w}]A_\ell\,[X^w]z_\ell.
\]
Thus there are an index $\ell$ and $w\in\supp z_\ell\subseteq V$ with $[X^{u-w}]A_\ell\ne0$. Comparing degree $u+v-c-w$ in $sA_\ell=rB_\ell$ gives
\[
[X^{v-w}]B_\ell=[X^{u-w}]A_\ell\ne0.
\]
Hence $u-w,v-w\in P$, that is, $w\le u,v$.
\end{proof}

\begin{cor}\label{cor:finite-global}
For $p,q\in P$, the $R$-module $J(p,q)$ is flat if and only if
\[
\forall\,\varnothing\ne F\subseteq U(p,q),\ |F|<\infty,\quad
\exists w\in U(p,q)\ \forall u\in F:\ w\le_Pu.
\]
Moreover, $J(p,q)$ is projective in $\Gr_LR$ if and only if $U(p,q)$ has a least element.
\end{cor}
\begin{proof}
By induction, downward directedness for pairs is equivalent to the stated finite-set condition, so the first assertion follows from Thm~\ref{thm:flat-ideal}. The second follows from Thm~\ref{thm:projective}.
\end{proof}

\subsection{Multiplicative obstructions in a general DPO}
\begin{prop}\label{prop:upper-multiplication}
For $a\in P\setminus\{0\}$,
\[
aU(p,q)+P\subseteq U(ap,aq).
\]
The left-hand side is the region obtained from the original common upper bounds by multiplication and positive translation; its ideal is $\phi_a(J(p,q))R$.
\end{prop}
\begin{proof}
If $v\in U(p,q)$ and $r\in P$, then
\[
av+r-ap=a(v-p)+r\in P,\qquad av+r-aq\in P.
\]
The monomial basis of $\phi_a(J(p,q))R$ has exponent set $aU(p,q)+P$, so this is the ideal of that region.
\end{proof}

\begin{thm}[General binary multiplication formula]\label{thm:general-binary}
For a DPO positive cone $P$, $p,q\in P$, and $a\in P\setminus\{0\}$,
\[
\Tor_1^R(R/(X^p,X^q),R_a)
\cong\frac{X^{ap}R\cap X^{aq}R}{\phi_a(J(p,q))R}.
\]
In the lattice case, this recovers Thm~\ref{thm:defect}.
\end{thm}
\begin{proof}
Let $F_1=R\sh p\oplus R\sh q$ and $d_1(f,g)=X^pf+X^qg$. There is an isomorphism
\[
J(p,q)\xrightarrow{\sim}\ker d_1,\qquad
h\mapsto(h/X^p,-h/X^q).
\]
Take the free module $F_2=\bigoplus_{w\in U(p,q)}R\sh w$. The map
\[
e_w\longmapsto(X^{w-p},-X^{w-q})
\]
surjects onto $\ker d_1$. Extend this initial segment to a free resolution. After tensoring with $R_a$,
\[
\ker d_1^a\xrightarrow{\sim}X^{ap}R\cap X^{aq}R,
\qquad(f,g)\mapsto X^{ap}f.
\]
This sends $d_2^a(e_w)$ to $X^{aw}$, so the image of the boundaries is
$\sum_{w\in U(p,q)}X^{aw}R=\phi_a(J(p,q))R$. Passing to the quotient gives the isomorphism.
\end{proof}

\begin{cor}\label{cor:upper-defect}
For $p,q\in P$ and $a\in P\setminus\{0\}$,
\[
\Tor_1^R(R/(X^p,X^q),R_a)=0
\iff U(ap,aq)=aU(p,q)+P.
\]
\end{cor}
\begin{proof}
The numerator and denominator have monomial $\Z$-bases indexed by $U(ap,aq)$ and $aU(p,q)+P$, respectively. Their quotient vanishes exactly when these sets coincide.
\end{proof}

\subsection{Topological criteria for lattice ordering}
\begin{prop}[Topological condition for doubling]\label{prop:double}
The map $m_2:X_P\to X_P$, $x\mapsto2x$, is a continuous bijection, and
\[
m_2\text{ is a self-homeomorphism}
\iff \bigl(\forall x\in L,\quad 2x\in P\Longrightarrow x\in P\bigr).
\]
\end{prop}
\begin{proof}
If $x\le y$, then $2(y-x)\in P$, so $m_2$ preserves order. Characteristic zero gives bijectivity.
A map between upper Alexandrov spaces is continuous if and only if it preserves order. Thus
\[
m_2^{-1}\text{ is continuous}\iff p/2\in P\ (\forall p\in P)
\iff(2x\in P\Rightarrow x\in P).
\]
\end{proof}

\begin{thm}[Complete topological criterion for lattice ordering]\label{thm:general-compact}
For any DPO positive cone $P$, the cone $P$ is a lattice if and only if
\begin{enumerate}
\item $m_2$ is a self-homeomorphism of $X_P$;
\item $U_x\cap U_y$ is compact for all $x,y\in L$.
\end{enumerate}
\end{thm}
\begin{proof}
If $P$ is a lattice, Prop~\ref{prop:lattice-basics} gives $P/2=P$, so $m_2$ is a self-homeomorphism. Also, $U_x\cap U_y=U_{x\vee y}$ is compact.
Conversely, (1) gives $P/2\subseteq P$. By (2) and Prop~\ref{prop:compact-cover}, each $U_x\cap U_y$ is a finite union of pairwise nonnested sets $U_{z_j}$. Choose the auxiliary total order constructed in the proof of Thm~\ref{thm:compact-min}. If there are at least two $z_j$, the midpoint of the two least ones remains in the intersection, since
\[
\tfrac12(z_1+z_2)-x\in P,\qquad\tfrac12(z_1+z_2)-y\in P.
\]
The same finite-minimal-element argument gives $z_1\le_Pz_2$, a contradiction. Thus the intersection is a single $U_z$, and Prop~\ref{prop:carrier-basics} gives a lattice order.
\end{proof}

\begin{exmp}[A space obtained by coordinate extension and rational layering]\label{ex:lattice-space}
Let $E=\Q$ with minimal open neighborhoods $N_E(r)=[r,\infty)\cap\Q$. Set
\[
B=\Q^{(\mathbb R)},\qquad
N_B(u)=\{v\in B:v_j\ge u_j\text{ for every }j\in\mathbb R\}.
\]
The elements of $B$ are finitely supported functions, $\mathbb R$ is only an index set, and all coordinate comparisons use the usual order on $\Q$.
For a $T_0$ Alexandrov space $Z$ with a distinguished zero, denote the minimal open neighborhood of $x$ by $N_Z(x)$ and define
\begin{align*}
\mathsf A(Z)&=Z\times B,&
N_{\mathsf A(Z)}(x,u)&=N_Z(x)\times N_B(u),\\
\mathsf L(Z)&=\Q\times Z,&
N_{\mathsf L(Z)}(r,x)&=\bigl(\{r\}\times N_Z(x)\bigr)
\cup\bigl((r,\infty)\cap\Q\bigr)\times Z.
\end{align*}
Embed the old space by $x\mapsto(x,0)$ and $x\mapsto(0,x)$, respectively. Each topology consists of arbitrary unions of the indicated minimal open neighborhoods. Give an increasing union $Y=\bigcup_\beta Z_\beta$ the topology
\[
O\subseteq Y\text{ is open}\iff O\cap Z_\beta\text{ is open in }Z_\beta\text{ for all }\beta.
\]

Fix an enumeration of $\Q$, and let $S(D)$ be the additive monoid generated by $D$. Starting from
\[
(X_0,a_0,D_0)=(E^3,1,\{1\})
\]
proceed recursively as follows, with $E^3$ carrying coordinatewise minimal open neighborhoods.
At a successor stage, set
\[
b_\alpha=\min\{s\in S(D_\alpha):s>a_\alpha\},
\]
and choose the first rational in the enumeration satisfying $a_\alpha<h_\alpha<b_\alpha$. Set
\[
(X_{\alpha+1},a_{\alpha+1},D_{\alpha+1})
=\bigl(\mathsf A(X_\alpha),h_\alpha,D_\alpha\cup\{h_\alpha\}\bigr).
\]
At a nonzero limit stage $\lambda$, set
\[
Y=\bigcup_{\beta<\lambda}X_\beta,\qquad
a=\sup_{\beta<\lambda}a_\beta,\qquad
D=\bigcup_{\beta<\lambda}D_\beta.
\]
If $a<\infty$, put $a_\lambda=a$ and $D_\lambda=D$, and set
\[
X_\lambda=
\begin{cases}
\mathsf L(Y),&a\in S(D),\\
Y,&a\notin S(D).
\end{cases}
\]
If $a=\infty$, stop and set $X=Y$.
The resulting $X$ is a connected $T_0$ Alexandrov space. Simultaneous doubling of all rational coordinates defines a self-homeomorphism $m_2$, and
\[
\forall x,y\in X\quad\exists!z\in X:\qquad N_X(x)\cap N_X(y)=N_X(z).
\]
It is homeomorphic to the order-carrying space obtained in \cite[\S5--6]{BY} on $\C$ from the three-dimensional initial configuration with the same choices of heights.
\end{exmp}
\begin{proof}
For well-ordering of $S(D)$ and termination of the recursion, see \byref{12}{Lem~4.3} and \byref{16}{Thm~5.4}.
At each recursive stage, we construct the unique point $x\sqcup y$ satisfying $N(x)\cap N(y)=N(x\sqcup y)$. On $E$ and $B$, it is given by the usual coordinatewise maximum. In the two extensions,
\[
(x,u)\sqcup(y,v)=(x\sqcup y,\max(u,v)),
\]
\[
(r,x)\sqcup(s,y)=
\begin{cases}
(r,x),&r>s,\\
(s,y),&s>r,\\
(r,x\sqcup y),&r=s.
\end{cases}
\]
The embeddings preserve these operations and the traces of minimal neighborhoods. In an increasing union, for $x\in X_\beta$,
\[
N_X(x)=\bigcup_{\gamma\ge\beta}N_{X_\gamma}(x).
\]
Thus the intersection formula persists, as do the $T_0$ and Alexandrov properties. Any two basic open sets intersect, giving connectedness. Doubling and division by $2$ are inverse homeomorphisms in each coordinate and layer, and commute with the embeddings.
Finally, the two extensions in \byref{9}{Lem~3.5} give the displayed neighborhood formulas. To check the number of coordinates in the Bao--Zhang construction, use its local notation: $T$ is a transcendence basis of $\C$, $t\in T$, $d$ is the height vanishing on $\Q(T\setminus\{t\})^\times$, $U\subseteq\{d\le A\}$ is the current subspace, and $L_{<h}=\{x\in\C:d(x)<h\}$. These data are given in \byref{10}{Lem~4.1} and \byref{14}{Prop~5.2}. For $r\in\Q\cap(A,h)$, the elements
\[
\{t^r\xi+U:\xi\in T\setminus\{t\}\}
\]
are linearly independent in $L_{<h}/U$, since every nonzero finite rational combination has height $r>A$. Thus the index set of new coordinates has cardinality $|T|=|\mathbb R|$ at every stage and can be identified with the coordinate set of $B$. Choose these labels stagewise and take unions at limits to obtain a bijection preserving all minimal neighborhoods. For coverage of the whole field, see \byref{16}{Thm~5.4} and \byref{18}{Thm~6.1}.
\end{proof}

\begin{cor}\label{cor:topological-obstruction}
If $m_2$ is a self-homeomorphism, then, for every $x\in L$,
\[
o_P(x)=0\iff V_x\text{ is compact}.
\]
\end{cor}
\begin{proof}
Since $m_2$ is a self-homeomorphism, the midpoint argument in the proof of Thm~\ref{thm:general-compact} applies separately to each $V_x=U_0\cap U_x$. Thus
\[
V_x\text{ is compact}\iff\min V_x\text{ exists}\iff o_P(x)=0,
\]
where the last step is Thm~\ref{thm:obstruction}.
\end{proof}

\section{Reconstruction of order-carrying spaces and the comparison theorem}
\subsection{Intrinsic characterization of order-carrying spaces}
\begin{defn}[Intrinsic order-carrying space]\label{def:intrinsic-carrier}
Fix the field $L$. A topology $\tau$ is order-carrying if:
\begin{enumerate}
\item arbitrary intersections of open sets are open;
\item $(L,\tau)$ is $T_0$;
\item each translation $x\mapsto x+c$ is a self-homeomorphism;
\item $(L,\tau)$ is connected;
\item field multiplication $L\times L\to L$ is continuous at $(0,0)$.
\end{enumerate}
The product carries the product topology. This definition does not require $L$ to be a topological field.
\end{defn}

\begin{thm}[Reconstruction of positive cones and order-carrying topologies]\label{thm:reconstruction}
On a fixed field $L$, DPO positive cones correspond bijectively to order-carrying topologies, with inverse maps
\[
P\longmapsto\tau_P,\qquad
\tau\longmapsto U_0,
\]
where $U_0$ is the minimal open neighborhood of zero. The correspondence retains the underlying field and its operations.
\end{thm}
\begin{proof}
By Prop~\ref{prop:carrier-basics}, $\tau_P$ satisfies the first four conditions. The set $P\times P$ is the minimal neighborhood of $(0,0)$, and $PP\subseteq P$, so multiplication is continuous there.
Conversely, set $P=U_0$. Translation homeomorphisms give $U_x=x+P$. If $p\in P$, then $U_p\subseteq P$, hence $P+P\subseteq P$. If $p,-p\in P$, then $0\in U_p$ and $p\in U_0$, so $U_p=U_0$; the $T_0$ property gives $p=0$.
The set $H=P-P$ is an additive subgroup, and, for each $c\in L$,
\[
c+H=\bigcup_{h\in H}U_{c+h}.
\]
Thus all cosets of $H$ are open. Since $H$ is nonempty and clopen, connectedness gives $H=L$. Continuity of multiplication at $(0,0)$ gives open neighborhoods $A,B$ of zero with $AB\subseteq P$; since $P\subseteq A\cap B$, we obtain $PP\subseteq P$.
Finally, $O\in\tau\iff\forall x\in O,\ x+P\subseteq O$, so $\tau=\tau_P$.
\end{proof}

\begin{exmp}[A space of strings with well-ordered support]\label{ex:dpo-space}
Let $A=\N_0\times\Q$, $o=(0,0)$, and take the point set
\[
X=\{f:\Q\to A:\ \supp f=\{r:f(r)\ne o\}\text{ is well ordered in the usual order}\}.
\]
Denote the zero string by $0$. For $f\ne g$, write $r(f,g)=\min\{r:f(r)\ne g(r)\}$. Prescribe minimal open neighborhoods
\[
N_X(f)=\{f\}\cup
\left\{g\in X\setminus\{f\}:
\begin{array}{l}
f(r(f,g))=(j,a),\quad g(r(f,g))=(j,b)\\
\text{for some }j\in\N_0\text{ and }a<b
\end{array}\right\},
\]
and set $\tau=\{\bigcup_{f\in C}N_X(f):C\subseteq X\}$.
This is a connected $T_0$ Alexandrov space. It admits field operations, isomorphic to those of $\C$, for which it is an order-carrying space. There are two minimal open neighborhoods whose intersection is not the minimal open neighborhood of any point.
\end{exmp}
\begin{proof}
The union of two well-ordered supports is well ordered, so $r(f,g)$ is defined. Comparing the first differing position gives
\[
g\in N_X(f)\Longrightarrow N_X(g)\subseteq N_X(f),\qquad
g\in N_X(f),\ f\in N_X(g)\Longrightarrow f=g.
\]
Hence the prescribed sets are minimal open neighborhoods and the space is $T_0$ Alexandrov. For any $f,g$, choose $s\in\Q$ below every element of both supports. The point with value $(0,1)$ at $s$ and $o$ elsewhere lies in $N_X(f)\cap N_X(g)$, proving connectedness.
Choose coset representatives $(c_j)_{j\in\N_0}$ for $\overline\Q/\Q$, with $c_0=0$. The map
\[
f(r)=(j,a)\longmapsto c_j+a
\]
identifies $X$ with the Hahn field $H=\overline\Q((t^\Q))$. The neighborhood condition is equivalent to $g=f$ or the leading coefficient of $g-f$ belonging to $\Q_{>0}$. Hence
\[
N_X(f+h)=h+N_X(f),\qquad N_X(0)N_X(0)\subseteq N_X(0).
\]
Thus translations are homeomorphisms and multiplication is continuous at $(0,0)$. By \cite{MacLane}, $H$ is algebraically closed. It has characteristic zero and cardinality $|\mathbb R|$, so $\operatorname{trdeg}_{\Q}H=|\mathbb R|$, and hence $H\cong\C$.
Let $\eta(0)=(1,0)$ and let $\eta$ equal $o$ elsewhere. The first-difference rule gives
\[
N_X(0)\cap N_X(\eta)
=\{g:\ s=\min\supp g<0,\ g(s)=(0,b),\ b>0\}.
\]
For any such $g$, choose $s<q<0$ and let $h$ have value $(0,1)$ at $q$ and $o$ elsewhere. Then $h$ remains in the intersection but $h\notin N_X(g)$, so the intersection is not the minimal open neighborhood of any point.
\end{proof}

\begin{cor}\label{cor:marked-field}
Let $P,Q$ be DPO positive cones on the fixed field $L$. Then
\[
\tau_P=\tau_Q\iff P=Q.
\]
\end{cor}
\begin{proof}
By Thm~\ref{thm:reconstruction}, $P$ and $Q$ are the minimal open neighborhoods of zero in their respective topologies, so equality of the topologies implies $P=Q$. The converse follows from the definition of the order-carrying topology.
\end{proof}

The underlying field and its operations are fixed here. Abstract homeomorphism types do not specify these data and are not the objects of this bijection.

\begin{cor}[Cone inclusion and refinement of topologies]\label{cor:cone-topology}
Let $P,Q$ be DPO positive cones on the same field $L$, and write $U_x^P=x+P$ and $U_x^Q=x+Q$. Then
\[
P\subseteq Q
\iff\tau_Q\subseteq\tau_P
\iff\id_L:X_P\longrightarrow X_Q\text{ is continuous}.
\]
Moreover, $P\subsetneq Q$ if and only if $\tau_Q\subsetneq\tau_P$.
\end{cor}
\begin{proof}
If $P\subseteq Q$ and $O\in\tau_Q$, then, for every $x\in O$,
$x+P\subseteq x+Q\subseteq O$, so $O\in\tau_P$.
Conversely, if $\tau_Q\subseteq\tau_P$, then $Q=U_0^Q$ is $\tau_P$-open and contains $0$, so $P=U_0^P\subseteq Q$.
The second equivalence is the definition of continuity. Applying the first equivalence in both directions gives the assertion about strict inclusion.
\end{proof}

\subsection{Truncated regions and singular homology}
\begin{defn}[Truncated region]\label{def:truncation}
Let $V\subseteq P$ be upward closed, $a\in P\setminus\{0\}$, and $u\in P$. Set
\[
Y_{a,u}(V)=\{v\in V:av\le_Pu\}\subseteq X_P,
\]
Give it the subspace topology, equivalently, the upper Alexandrov topology of its induced order.
\end{defn}

\begin{prop}[Order-chain resolution of a region ideal]\label{prop:bar}
Set
\[
B_k(V)=\bigoplus_{v_0<_P\cdots<_Pv_k\text{ in }V}R\sh{v_k}.
\]
For $k\ge1$, define
\begin{align*}
d[v_0,\ldots,v_k]
&=\sum_{i=0}^{k-1}(-1)^i[v_0,\ldots,\widehat v_i,\ldots,v_k]\\
&\quad+(-1)^kX^{v_k-v_{k-1}}[v_0,\ldots,v_{k-1}],
\end{align*}
and let $B_0(V)\to J_V$ be given by $[v]\mapsto X^v$. This is a graded free resolution of $J_V$.
\end{prop}
\begin{proof}
For each $w\in L$, the homogeneous basis associated to the chain labels is
\[
X^{w-v_k}[v_0,\ldots,v_k],\qquad v_k\le_P w.
\]
In this basis the differential is the usual face-deletion boundary, so $d^2=0$. If $w\notin V$, all terms vanish. If $w\in V$, the augmented degree component identifies with
\[
C_\bullet(\ord(V\cap\downarrow w);\Z)\longrightarrow\Z.
\]
where $\downarrow w=\{v:v\le_Pw\}$, and $V\cap\downarrow w$ has greatest element $w$. The vertex-adjoining operator in the proof of Prop~\ref{prop:directed-acyclic} gives a contraction, proving exactness in every degree.
\end{proof}

\begin{thm}[Order-chain comparison for a general DPO]\label{thm:tor-order}
Let $P$ be a DPO positive cone, $V\subseteq P$ upward closed, $a\in P\setminus\{0\}$, and $u\in P$. For $n\ge1$, there are natural isomorphisms
\[
\Tor_n^R(R/J_V,R_a)_u
\cong\widetilde H_{n-1}(\ord(Y_{a,u}(V));\Z).
\]
\end{thm}
\begin{proof}
Splicing $B_\bullet(V)\to J_V$ with $J_V\hookrightarrow R$ gives a free resolution of $R/J_V$ with $F_0=R$ and $F_n=B_{n-1}(V)$ for $n\ge1$. For $n\ge1$, define the basis bijection
\[
[v_0,\ldots,v_{n-1}]\otimes X^{u-av_{n-1}}
\longmapsto[v_0,\ldots,v_{n-1}].
\]
The left-hand side exists if and only if $av_{n-1}\le_Pu$. Since $aP\subseteq P$, this is equivalent to all $v_j$ belonging to $Y_{a,u}(V)$. The bijection sends $d^a$ to the face-deletion boundary and $F_1\to F_0$ to the sum of coefficients; $u\in P$ ensures $(F_0\otimes_RR_a)_u=\Z$. Thus
\[
(F_\bullet\otimes_RR_a)_u\cong
\widetilde C_{\bullet-1}(\ord(Y_{a,u}(V));\Z),
\]
Passing to homology gives the result. Chain labels are unchanged, giving naturality.
\end{proof}

\begin{prop}[Order complexes and Alexandrov spaces]\label{prop:mccord}
Give a poset $A$ the upper Alexandrov topology with minimal open neighborhoods $\up x=\{y:x\le y\}$. There is a natural weak homotopy equivalence
\[
\mu_A:|\ord(A)|\longrightarrow A,
\]
sending each point to the greatest vertex of its support simplex. It induces isomorphisms on all singular homology groups. No finiteness or countability assumption on $A$ is required.
\end{prop}
\begin{proof}
Let $\up x=\{y:x\le y\}$, and write $(\lambda_v)$ for barycentric coordinates in the realization. Then
\[
\mu_A^{-1}(\up x)=\left\{\lambda:\sum_{v\ge x}\lambda_v>0\right\}
\]
is open in each simplex and hence in the realization. Set $s(\lambda)=\sum_{v\ge x}\lambda_v$ and, on this inverse image, define
\[
r_x(\lambda)_v=
\begin{cases}\lambda_v/s(\lambda),&v\ge x,\\0,&v\not\ge x.\end{cases}
\]
The map $H_t(\lambda)=(1-t)\lambda+t r_x(\lambda)$ is a deformation retraction onto $|\ord(\up x)|$, remaining in the original support simplex at every stage. The latter is a cone with apex $x$ and hence contractible; $\up x$ is also contractible because it has a least element.
Thus the restriction over each minimal open set is a weak homotopy equivalence. McCord's basis criterion \cite[Theorem 6]{McCord} implies that $\mu_A$ is a weak homotopy equivalence. For an order-preserving map $f:A\to B$,
\[
\mu_B\circ|\ord(f)|=f\circ\mu_A,
\]
since $f(\max\supp\lambda)=\max f(\supp\lambda)$. Weak homotopy equivalences induce isomorphisms on singular homology.
\end{proof}

\begin{thm}[Realization by ordinary singular homology]\label{thm:tor-singular}
Let $P$ be a DPO positive cone, $V\subseteq P$ upward closed, $a\in P\setminus\{0\}$, and $u\in P$. For $n\ge1$, there are natural isomorphisms
\[
\Tor_n^R(R/J_V,R_a)_u
\cong\widetilde H_{n-1}(Y_{a,u}(V);\Z).
\]
For $r\in P$, the action of $X^r$ corresponds to $Y_{a,u}(V)\hookrightarrow Y_{a,u+r}(V)$. For $V\subseteq W$, the quotient map $R/J_V\to R/J_W$ corresponds to inclusion of truncated regions.
\end{thm}
\begin{proof}
Compose the natural isomorphisms
\[
\Tor_n^R(R/J_V,R_a)_u
\xrightarrow{\sim}\widetilde H_{n-1}(\ord(Y_{a,u}(V));\Z)
\xrightarrow{\sim}\widetilde H_{n-1}(Y_{a,u}(V);\Z),
\]
where the second step uses simplicial--singular comparison and Prop~\ref{prop:mccord}. Multiplication by $X^r$ sends $[v_0,\ldots,v_k]\otimes X^{u-av_k}$ to the basis element with the same label in degree $u+r$. The chain map for $V\subseteq W$ also preserves labels. These therefore correspond to the stated inclusions.
\end{proof}

\begin{cor}\label{cor:global-singular}
For every nonempty open set $O\subseteq X_P$ and $n\ge0$,
\[
\widetilde H_n(O;\Z)=0.
\]
In particular, $\widetilde H_n(X_P;\Z)=0$ for every $n\ge0$.
\end{cor}
\begin{proof}
A nonempty open set is upward directed, so Prop~\ref{prop:directed-acyclic} and Prop~\ref{prop:mccord} give
\[
\widetilde H_n(O;\Z)\cong\widetilde H_n(\ord(O);\Z)=0\quad(n\ge0).
\]
Take $O=X_P$ for the assertion about the whole space.
\end{proof}

Truncated regions need not be open, so this vanishing result does not apply to them directly.

\begin{cor}[Homological comparison for extensions of positive cones]\label{cor:extension-comparison}
Let $P\subseteq Q$ be DPO positive cones on the same field $L$. Write $R_P=\Z[P]$ and $R_Q=\Z[Q]$, and give $R_Q$ its module structure through $R_P\hookrightarrow R_Q$.
For a $P$-upward-closed subset $V\subseteq P$ and $u\in Q$, set
\[
J_V^P=\bigoplus_{v\in V}\Z X^v,\qquad
Y_{P,Q,u}(V)=V\cap(u-Q)=V\cap\overline{\{u\}}^{\,X_Q}.
\]
Give $Y_{P,Q,u}(V)$ the subspace topology of $X_P$ and form its order complex using the induced order $\le_P$. For $n\ge1$, there are natural isomorphisms
\begin{align*}
\Tor_n^{R_P}(R_P/J_V^P,R_Q)_u
&\cong\widetilde H_{n-1}(\ord(Y_{P,Q,u}(V));\Z)\\
&\cong\widetilde H_{n-1}(Y_{P,Q,u}(V);\Z).
\end{align*}
Both rings are $L$-graded. For $r\in Q$, right multiplication by $X^r$ corresponds to $Y_{P,Q,u}(V)\hookrightarrow Y_{P,Q,u+r}(V)$. For $V\subseteq W$, the quotient map corresponds to inclusion of truncated regions.
\end{cor}
\begin{proof}
By Prop~\ref{prop:bar}, the modules $F_0=R_P$ and $F_n=B_{n-1}^P(V)$ for $n\ge1$ form a free resolution of $R_P/J_V^P$, where the superscript $P$ denotes the order-chain resolution for the cone $P$.
After tensoring with $R_Q$, the degree-$u$ component of the $n$th term has basis
\[
[v_0,\ldots,v_{n-1}]\otimes X^{u-v_{n-1}},\qquad
v_0<_P\cdots<_Pv_{n-1},\quad u-v_{n-1}\in Q.
\]
Since $P\subseteq Q$, the final condition is equivalent to all vertices belonging to $Y_{P,Q,u}(V)$.
Send this basis element to $[v_0,\ldots,v_{n-1}]$. The differential becomes the face-deletion boundary, the map from term $1$ to term $0$ becomes the augmentation, and $(R_Q)_u=\Z$. Passing to homology gives the first isomorphism; Prop~\ref{prop:mccord} gives the second.
Multiplication by $X^r$ and the chain map for $V\subseteq W$ preserve vertex labels and hence correspond to the stated inclusions. This naturality is preserved under Prop~\ref{prop:mccord}.
\end{proof}

\subsection{Finite nerve models for lattices}
\begin{thm}[Nerve model of a truncated region]\label{thm:nerve}
Let $P$ be a lattice, $m\ge1$, $\mathbf p=(p_1,\ldots,p_m)\in P^m$, $a\in P\setminus\{0\}$, and $u\in P$. Set $V=\bigcup_{i=1}^m(p_i+P)$ and write
\[
Y_i=\{v:p_i\le_Pv,\ av\le_Pu\}.
\]
Then
\[
|\ord(Y_{a,u}(V))|\simeq|\Xi_{a,u}(\mathbf p)|.
\]
The natural map $\mu_Y:|\ord(Y)|\to Y$, with $Y=Y_{a,u}(V)$, is a weak homotopy equivalence. Consequently, $Y$ and $|\Xi_{a,u}(\mathbf p)|$ have isomorphic singular homology groups.
\end{thm}
\begin{proof}
Let $Y=Y_{a,u}(V)$. The least vertex of each nonempty chain in $Y$ belongs to some $Y_i$. Since $Y_i$ is upward closed in $Y$, the whole chain lies in $Y_i$. Thus $\{\ord(Y_i)\}_{i=1}^m$ is a subcomplex cover.
For $\varnothing\ne I\subseteq[m]$,
\[
\bigcap_{i\in I}Y_i=\{v:b_I\le_Pv,\ av\le_Pu\}\ne\varnothing
\iff ab_I\le_Pu.
\]
When nonempty, the intersection has least element $b_I$, so its order complex is contractible. By the nerve theorem for finite subcomplex covers \cite[Theorem 10.6]{Bjorner}, its nerve is precisely $\Xi_{a,u}(\mathbf p)$, giving the asserted homotopy equivalence. If $Y=\varnothing$, both realizations are empty. The singular homology assertion follows from Prop~\ref{prop:mccord}.
\end{proof}

\begin{rem}\label{rem:finite-model}
For a lattice $P$ and $V=\bigcup_{i=1}^m(p_i+P)$, Thm~\ref{thm:nerve} reduces the homology of the truncated region to that of the finite complex $\Xi_{a,u}(\mathbf p)$, whose nonempty faces are determined by $ab_I\le_Pu$. For a general DPO positive cone, the full order-complex model of Thm~\ref{thm:tor-order} remains applicable.
\end{rem}

\subsection{The unified comparison theorem and its naturality}
\begin{thm}[Unified comparison theorem]\label{thm:main}
Fix a field $L$ of characteristic zero. DPO positive cones and order-carrying topologies determine one another. For the corresponding $P$, $X_P$, and $R=\Z[P]$, we have:
\begin{enumerate}
\item an exact equivalence of categories $\Gr_LR\simeq\Sh(X_P;\Ab)$;
\item for upward-closed $V\subseteq P$ and $M\in\Gr_LR$,
\[
\Ext^n_{\Gr_LR}(J_V,M)\cong H^n(V;\widetilde M|_V),\qquad n\ge0;
\]
\item for upward-closed $V\subseteq P$, $a\in P\setminus\{0\}$, and $u\in P$,
\[
\Tor_n^R(R/J_V,R_a)_u
\cong\widetilde H_{n-1}(Y_{a,u}(V);\Z),\qquad n\ge1.
\]
\end{enumerate}
These correspondences are compatible with degree-preserving module maps, inclusions of regions, positive translations of the upper bound, and field isomorphisms preserving positive cones. Finitely generated regions in the lattice case also admit the finite nerve models of Thm~\ref{thm:nerve}.
\end{thm}
\begin{proof}
The object correspondence is Thm~\ref{thm:reconstruction}; (1)--(2) are Thm~\ref{thm:sheaf-equivalence}; and (3) is Thm~\ref{thm:tor-singular}.
If $\sigma:L\to L'$ is a field isomorphism with $\sigma(P)=P'$, then
\[
X^v\mapsto X^{\sigma(v)},\qquad
[v_0,\ldots,v_k]\mapsto[\sigma(v_0),\ldots,\sigma(v_k)]
\]
give isomorphisms of rings and free resolutions, respectively, and send $Y_{a,u}(V)$ to $Y_{\sigma(a),\sigma(u)}(\sigma(V))$. They commute with differentials, restriction maps, and McCord maps. The remaining naturality follows from the explicit constructions in Thm~\ref{thm:sheaf-equivalence} and Thm~\ref{thm:tor-singular}; finite compression follows from Thm~\ref{thm:nerve}.
\end{proof}

\par\medskip
\noindent\textbf{Important Rmk.} After completing the initial version of this manuscript, the author became aware, through a subsequent literature search, of the following general result, recorded in \cite[Theorem~2.21]{BM}:

Let $G$ be a partially ordered abelian group, let $G_+=\{g\in G:0\le g\}$, and let $A$ be a commutative ring with identity. Equip $G$ with its upper Alexandrov topology, denoted by $G_{\mathrm{Alex}}$. There are natural exact equivalences
\[
\Gr_G A[G_+]
\simeq
\operatorname{Fun}(G,A\text{-}\mathrm{Mod})
\simeq
\Sh(G_{\mathrm{Alex}};A\text{-}\mathrm{Mod}),
\]
where graded morphisms preserve degree, and the stalk of the sheaf corresponding to a graded module $M$ at $g$ is $M_g$.

The constructions already established in this paper supply the following data:
\begin{enumerate}
\item The ordered additive group $G=(L,+,\le_P)$, its positive cone $G_+=P$, and the graded monoid ring $R=\Z[P]$ (Defs~\ref{def:dpo}--\ref{def:ring}).
\item The upper Alexandrov space $X_P$, with minimal neighborhoods $U_x=x+P$, and its intrinsic reconstruction from the topology (Def~\ref{def:carrier} and Thm~\ref{thm:reconstruction}).
\item The region ideals $J_V=\bigoplus_{v\in V}\Z X^v$ (Def~\ref{def:monomial}), corresponding to the constant $\Z$-valued diagram on $V$, extended by zero.
\item The coefficient modules $R_a$ associated with positive multiplication and the truncated regions $Y_{a,u}(V)=\{v\in V:av\le_Pu\}$ (Defs~\ref{def:twist} and~\ref{def:truncation}).
\end{enumerate}

Substituting $A=\Z$ and these data into the general equivalence immediately gives Thm~\ref{thm:main}(1). The identification
\[
\Hom_{\Gr_LR}(J_V,M)
\cong \Gamma(V,\widetilde M|_V)
\]
then yields part~(2) by passage to derived functors. For part~(3), the change of coefficients $M\mapsto M\otimes_RR_a$ corresponds to left Kan extension along $f_a(v)=av$, whose indexing poset over $u$, restricted to $V$, is precisely $Y_{a,u}(V)$. The standard formula for derived Kan extensions \cite[\S1.3]{Dj}, together with
\[
0\longrightarrow J_V\longrightarrow R\longrightarrow R/J_V\longrightarrow0
\]
and the McCord comparison in Prop~\ref{prop:mccord}, gives part~(3). Together with Thms~\ref{thm:reconstruction} and~\ref{thm:nerve}, this provides an abstract route to the unified comparison theorem.

After considerable deliberation, the author has chosen to retain both this abstract perspective and the explicit development above. The preceding proof follows, step by step, the theory developed in this paper; its constructions and chain-level identifications are fully explicit and closely tailored to the particular setting of DPOs.
\par\medskip

\begin{cor}[Three formulations of lattice ordering]\label{cor:three-lattice}
For each $x\in L$, the conditions $o_P(x)=0$, graded projectivity of $J_x$, and existence of a least element in $V_x$ are equivalent. They hold for all $x$ if and only if $P$ is a lattice; in this case, the order-carrying space also satisfies the complete topological criterion of Thm~\ref{thm:general-compact}.
\end{cor}
\begin{proof}
For each $x$, Thm~\ref{thm:obstruction} and Thm~\ref{thm:projective} give
\[
o_P(x)=0\iff\min V_x\text{ exists}\iff J_x\text{ is graded projective}.
\]
Imposing this condition for all $x$ gives a lattice by Prop~\ref{prop:reduced-min}; its topological characterization is Thm~\ref{thm:general-compact}.
\end{proof}

The next four corollaries concern changes of positive cone. For $P\subseteq Q$ on the same field, continue to write $R_P=\Z[P]$ and $R_Q=\Z[Q]$.

\begin{cor}[Nonflatness of strict extensions]\label{cor:extension-nonflat}
Let $P\subseteq Q$ be DPO positive cones. Then
\[
R_Q\text{ is flat over }R_P\iff P=Q.
\]
More precisely, for any $a\in Q\setminus P$ and representation $a=p-q$ with $p,q\in P$,
\[
\Tor_1^{R_P}\bigl(R_P/(X^p,X^q),R_Q\bigr)_p\cong\Z.
\]
\end{cor}
\begin{proof}
Such a representation exists because $P-P=L$. Set
\[
Y=\bigl((p+P)\cup(q+P)\bigr)\cap(p-Q),\qquad
Y_p=(p+P)\cap(p-Q),\quad Y_q=(q+P)\cap(p-Q).
\]
If $v\in Y_p$, then $v-p\in P\subseteq Q$ and $p-v\in Q$. Pointedness gives $v=p$, so $Y_p=\{p\}$.
Also, $q\in Y_q$ is its $\le_P$-least element, and $p\notin Y_q$.
For $v\in Y_q$, the inequality $v\le_Pp$ would imply $a=(p-v)+(v-q)\in P$, while $p\le_Pv$ would imply $v\in Y_p$. Both are impossible.
Thus $\ord(Y)=\ord(\{p\})\sqcup\ord(Y_q)$ is the disjoint union of two nonempty contractible complexes. Cor~\ref{cor:extension-comparison} gives
\[
\Tor_1^{R_P}\bigl(R_P/(X^p,X^q),R_Q\bigr)_p
\cong\widetilde H_0(\ord(Y);\Z)\cong\Z.
\]
Strict inclusion is therefore nonflat. If $P=Q$, then $R_Q$ is free of rank one over $R_P$.
\end{proof}

\begin{cor}[Extension criterion for preserving a reduced decomposition]\label{cor:extension-reduced}
Let $P\subseteq Q$ be DPO positive cones, and let $p,q\in P$ satisfy $\HH_P(p,q)=0$. There is a natural isomorphism of $L$-graded $R_Q$-modules
\[
\HH_Q(p,q)\cong
\Tor_1^{R_P}\bigl(R_P/(X^p,X^q),R_Q\bigr).
\]
Consequently, the following are equivalent:
\begin{enumerate}
\item $p-q$ remains a reduced decomposition in $Q$;
\item $\Tor_1^{R_P}(R_P/(X^p,X^q),R_Q)=0$;
\item $U_p^Q\cap U_q^Q=U_{p+q}^Q$.
\end{enumerate}
\end{cor}
\begin{proof}
By Thm~\ref{thm:binary} and injectivity of monomial multiplication, the following is a graded free resolution:
\[
0\longrightarrow R_P\sh{p+q}
\xrightarrow{\,(-X^q,X^p)\,}
R_P\sh p\oplus R_P\sh q
\xrightarrow{\,(X^p,X^q)\,}R_P
\longrightarrow R_P/(X^p,X^q)\longrightarrow0.
\]
Tensoring its free terms with $R_Q$ gives the binary complex $K_Q(p,q)$. Its first homology is therefore $\HH_Q(p,q)$ and computes the stated $\Tor_1$.
The criterion of Thm~\ref{thm:binary} for $Q$ and its monomial description give
\[
\HH_Q(p,q)=0
\iff(p+Q)\cap(q+Q)=p+q+Q
\iff p-q\text{ is reduced in }Q,
\]
proving the equivalence of the three conditions.
\end{proof}

\begin{cor}[Strict extensions of a lattice destroy reduced decompositions]\label{cor:lattice-extension}
Let $P$ be a lattice on $L$ and $Q$ a DPO positive cone on the same field, with $P\subsetneq Q$.
For every $a\in Q\setminus P$, its reduced decomposition $a=p-q$ in $P$ satisfies
\[
\HH_P(p,q)=0,\qquad\HH_Q(p,q)_p\cong\Z.
\]
In particular, $Q$ does not preserve all reduced decompositions of $P$.
\end{cor}
\begin{proof}
The reduced decomposition exists by Prop~\ref{prop:reduced-min}.
Applying Cor~\ref{cor:extension-reduced} and then Cor~\ref{cor:extension-nonflat} gives
\[
\HH_Q(p,q)_p
\cong\Tor_1^{R_P}\bigl(R_P/(X^p,X^q),R_Q\bigr)_p
\cong\Z.
\]
Thm~\ref{thm:binary} then excludes reducedness in $Q$.
\end{proof}

\begin{cor}[Directed limits and vanishing of reducedness obstructions]\label{cor:cone-limit}
Let $I$ be a nonempty directed poset, and let $(P_\alpha)_{\alpha\in I}$ be DPO positive cones on the same field $L$, with $\alpha\le\beta$ implying $P_\alpha\subseteq P_\beta$. Set $P_\infty=\bigcup_\alpha P_\alpha$. Then $P_\infty$ is a DPO positive cone, and
\[
\tau_{P_\infty}=\bigcap_{\alpha\in I}\tau_{P_\alpha}.
\]
Fix $\alpha_0\in I$ and $p,q\in P_{\alpha_0}$. The ring inclusions induce a natural isomorphism of $L$-graded abelian groups
\[
\varinjlim_{\alpha\ge\alpha_0}\HH_{P_\alpha}(p,q)
\xrightarrow{\sim}\HH_{P_\infty}(p,q).
\]
Thus $p-q$ is reduced in $P_\infty$ if and only if
\[
\forall\alpha\ge\alpha_0\ \forall\xi\in\HH_{P_\alpha}(p,q)\quad
\exists\beta\ge\alpha:\quad
\xi\longmapsto0\text{ in }\HH_{P_\beta}(p,q).
\]
\end{cor}
\begin{proof}
Every finite subset of $P_\infty$ lies in some $P_\alpha$, so additive and multiplicative closure and pointedness are inherited from that stage. The equality $P_{\alpha_0}-P_{\alpha_0}=L$ gives $P_\infty-P_\infty=L$.
For $O\subseteq L$,
\[
O\in\tau_{P_\infty}
\iff\forall x\in O,\ x+P_\infty\subseteq O
\iff\forall\alpha,\ O\in\tau_{P_\alpha}.
\]
Write $K_\alpha=K_{P_\alpha}(p,q)$ and $K_\infty=K_{P_\infty}(p,q)$. Since each ring element has only finitely many monomials,
\[
R_{P_\infty}=\varinjlim_{\alpha\ge\alpha_0}R_{P_\alpha},\qquad
K_\infty=\varinjlim_{\alpha\ge\alpha_0}K_\alpha
\]
holds termwise, and all termwise maps are injective. Finite support gives
\[
Z_1(K_\infty)=\bigcup_{\alpha\ge\alpha_0}Z_1(K_\alpha),\qquad
B_1(K_\infty)=\bigcup_{\alpha\ge\alpha_0}B_1(K_\alpha),
\]
where $Z_1=\ker d_1$ and $B_1=\operatorname{im}d_2$. Taking quotients gives the homology isomorphism. For $z\in Z_1(K_\alpha)$,
\[
[z]=0\text{ in }H_1(K_\infty)
\iff\exists\beta\ge\alpha:\ z\in B_1(K_\beta).
\]
Combining this with Thm~\ref{thm:binary} gives the final equivalence.
\end{proof}

\subsection{Positive multiplication: flatness and positive inverses}
We now give a complete criterion for flatness of positive-multiplication coefficient modules. For a lattice, it also characterizes the multipliers preserving all reduced decompositions.

\begin{thm}[Flatness, positive inverses, and self-homeomorphisms]\label{thm:flat-mult}
Let $P$ be any DPO positive cone and $a\in P\setminus\{0\}$. Then
\[
R_a\text{ is flat}\iff a^{-1}\in P\iff m_a:X_P\to X_P\text{ is a self-homeomorphism}.
\]
If $P$ is a lattice, these conditions are also equivalent to $m_a$ preserving all lattice operations, and to $ap-aq$ remaining reduced whenever $p-q$ is reduced, for all $p,q\in P$.
\end{thm}
\begin{proof}
If $a^{-1}\in P$, then $aP=P$ and $\phi_a$ is a ring automorphism. The map $R_a\to R$, $f\mapsto\phi_a^{-1}(f)$, is a left-module isomorphism, so $R_a$ is free. Both $m_a$ and its inverse preserve order, so $m_a$ is a self-homeomorphism.
If $m_a$ is a self-homeomorphism, its inverse preserves order, so $a^{-1}P\subseteq P$. Since $a\in P$, we obtain $1\in P$ and then $a^{-1}\in P$.
Suppose $a^{-1}\notin P$. There is $b\in P$ with $a^{-1}b\notin P$; otherwise $a^{-1}P\subseteq P$ would give a contradiction. Write $a^{-1}b=p-q$ with $p,q\in P$, and take
\[
V=(p+P)\cup(q+P),\qquad u=ap,\qquad Y=Y_{a,u}(V).
\]
Let $Y_p=Y\cap(p+P)$ and $Y_q=Y\cap(q+P)$. If $v\in Y_p$, then $a(v-p)\in P\cap(-P)$, so $v=p$. Also, $q\in Y_q$ is its least element.
For $v\in Y_q$, the inequality $v\le p$ would give $p-q\in P$, a contradiction. If $p\le v$, then $v\in Y_p$ gives $v=p$, again a contradiction. Thus $\ord(Y)$ is the disjoint union of two nonempty contractible complexes, and
\[
\Tor_1^R(R/(X^p,X^q),R_a)_{ap}
\cong\widetilde H_0(Y;\Z)\cong\Z.
\]
Hence $R_a$ is not flat.
In the lattice case, if the bijection $m_a$ preserves lattice operations, then
$m_a(x)\le m_a(y)\Rightarrow m_a(x\wedge y)=m_a(x)\Rightarrow x\le y$, so its inverse preserves order. The converse follows because order automorphisms preserve suprema and infima.
Preservation of lattice operations implies preservation of reduced decompositions. Conversely, if all reduced decompositions are preserved, their uniqueness gives, for any $x,y\in L$,
\[
(ax)^+=ax^+,\qquad
a(x\vee y)=ay+a(x-y)^+=ay+(ax-ay)^+=ax\vee ay.
\]
The identity $x\wedge y=x+y-(x\vee y)$ then shows that $m_a$ preserves infima.
\end{proof}

\begin{cor}[Completeness of binary tests]\label{cor:flat-detection}
The module $R_a$ is flat if and only if $\Tor_1^R(R/(X^p,X^q),R_a)=0$ for all $p,q\in P$. In the lattice case, this is also equivalent to $\delta_a(p,q)=0$ for all $p,q\in P$.
\end{cor}
\begin{proof}
Flatness annihilates all positive-degree $\Tor$. The converse follows from the nonzero binary class in the proof of Thm~\ref{thm:flat-mult}. In the lattice case, Thm~\ref{thm:defect} gives
\[
\Tor_1^R(R/(X^p,X^q),R_a)=0
\iff R/(X^{\delta_a(p,q)})=0
\iff\delta_a(p,q)=0,
\]
The last step uses $X^dR=R\iff d=0$.
\end{proof}

\part{Structural criteria and applications}
\section{Structural applications of multiplicative obstructions}
\subsection{Multiplicative obstructions induced by the imaginary unit}
In this subsection, $L$ contains an element satisfying $i^2=-1$; in particular, the results apply to $L=\C$.

\begin{prop}\label{prop:i-incomparable}
Every DPO positive cone $P$ satisfies $i,-i\notin P$. If $P$ is a lattice, the reduced decomposition $i=p-q$ satisfies $p,q>_P0$.
\end{prop}
\begin{proof}
If $\eps i\in P$ for $\eps\in\{1,-1\}$, then $-1=(\eps i)^2\in P$ and $1=(-1)^2\in P$, contradicting pointedness. If $q=0$ or $p=0$ in the reduced decomposition, then $i\in P$ or $-i\in P$, respectively. Hence $p,q\ne0$.
\end{proof}

\begin{thm}[A universal nonzero multiplicative class]\label{thm:universal-class}
Let $P$ be a lattice and $i=p-q$ a reduced decomposition. Set
\[
a=p+q,\qquad u=p^2+pq+q^2,\qquad J=(X^p,X^q).
\]
Then
\[
\Tor_1^R(R/J,R_a)_u\cong\Z.
\]
Under the comparison isomorphism, the reduced zero-dimensional class $[p]-[q]$ corresponds to a generator of this group. Multiplication by $X^{pq}$ sends it to zero in degree $u+pq$.
\end{thm}
\begin{proof}
By Prop~\ref{prop:i-incomparable}, $p,q,pq>_P0$. Reducedness gives $p\vee q=p+q$. We have
\[
u-ap=q^2\in P,\qquad u-aq=p^2\in P,\qquad
u-a(p\vee q)=-pq\notin P.
\]
Thus $\Xi_{a,u}(p,q)$ consists of two discrete vertices, and Thm~\ref{thm:xi} gives $\Tor_1\cong\widetilde H_0(S^0;\Z)=\Z$.
The two pieces of the truncated region have least elements $p,q$, respectively, and are disjoint, so $[p]-[q]$ is a generator.
Since $u+pq=a(p\vee q)$, the complex $\Xi_{a,u+pq}$ is a full edge and has zero reduced homology. By naturality in Thm~\ref{thm:xi}, $X^{pq}$ annihilates the class.
\end{proof}

\begin{cor}\label{cor:universal-nonflat}
Suppose $L$ contains an element satisfying $i^2=-1$, and let $P$ be a lattice on $L$. Then there is $a\in P\setminus\{0\}$ such that $a^{-1}\notin P$ and $m_a$ does not preserve all lattice operations.
\end{cor}
\begin{proof}
Thm~\ref{thm:universal-class} gives $\Tor_1^R(R/J,R_a)\ne0$, so $R_a$ is not flat. By Thm~\ref{thm:flat-mult}, $a^{-1}\notin P$ and $m_a$ does not preserve all lattice operations.
\end{proof}

\subsection{Flatness and reducedness obstructions for inverse-closed cones}
\begin{defn}[Inverse-closed positive cone]\label{def:inverse-closed}
A cone $P$ is inverse closed if $a\in P\setminus\{0\}$ always implies $a^{-1}\in P$. For any DPO positive cone $P_0$, define
\[
P_0^{\mathrm{inv}}=\{p/q:p\in P_0,\ q\in P_0\setminus\{0\}\}.
\]
\end{defn}

\begin{prop}\label{prop:inverse-extension}
The set $P_0^{\mathrm{inv}}$ is an inverse-closed DPO positive cone containing $P_0$.
\end{prop}
\begin{proof}
We have
\[
\frac pq+\frac rs=\frac{ps+rq}{qs},\qquad
\frac pq\frac rs=\frac{pr}{qs},\qquad
\left(\frac pq\right)^{-1}=\frac qp\quad(p\ne0),
\]
Thus the set is closed under addition, multiplication, and inversion of nonzero elements. Choose $c\in P_0\setminus\{0\}$. Since $p=(pc)/c$, we have $P_0\subseteq P_0^{\mathrm{inv}}$ and hence $P_0^{\mathrm{inv}}-P_0^{\mathrm{inv}}=L$.
If $p/q=-r/s$, then $ps=-rq\in P_0\cap(-P_0)$, so $ps=rq=0$. Nonzero denominators give $p=r=0$, proving pointedness.
\end{proof}

\begin{thm}[Common-upper-bound modules of inverse-closed cones]\label{thm:inverse-interpolation}
If $P$ is inverse closed, every common-upper-bound module $J(p,q)$ is flat. If $L$ contains $i$, then $P$ is not a lattice and $o_P(i)\ne0$.
\end{thm}
\begin{proof}
Take $u,v\in U(p,q)$ and let $d=u+v-p-q\in P$. If $d=0$, then $(u-p)+(v-q)=0$, so pointedness gives $u=p,v=q$. Together with $q\le u,p\le v$, this makes all four points equal.
If $d\ne0$, set $w=(uv-pq)/d$. Since $d^{-1}\in P$, the identities
\[
\begin{aligned}
w-p&=(u-p)(v-p)/d,&w-q&=(u-q)(v-q)/d,\\
u-w&=(u-p)(u-q)/d,&v-w&=(v-p)(v-q)/d
\end{aligned}
\]
give $p,q\le w\le u,v$. Thus $U(p,q)$ is downward directed, and $J(p,q)$ is flat by Thm~\ref{thm:flat-ideal}.
If $i=r-s$ is reduced, Prop~\ref{prop:i-incomparable} gives $r,s>_P0$. Set $h=rs/(r+s)>_P0$. Then
\[
r-h=r^2/(r+s)\in P,\qquad s-h=s^2/(r+s)\in P.
\]
Thus $i=(r-h)-(s-h)$ is a positive representation. Reducedness would require $(r-h)-r=-h\in P$, contradicting pointedness. Therefore $i$ has no reduced decomposition, and Thm~\ref{thm:obstruction} gives $o_P(i)\ne0$.
\end{proof}

\begin{cor}\label{cor:tor-not-lattice}
Let $P$ be an inverse-closed DPO positive cone. Then $R_a$ is flat for every $a\in P\setminus\{0\}$, and, for every $R$-module $M$ and $n\ge1$,
\[
\Tor_n^R(M,R_a)=0.
\]
If $i\in L$ satisfies $i^2=-1$, then $o_P(i)\ne0$.
\end{cor}
\begin{proof}
Inverse closure and Thm~\ref{thm:flat-mult} give
\[
\Tor_n^R(M,R_a)=0\quad(n\ge1)
\]
for all $M$ and $a>_P0$, while Thm~\ref{thm:inverse-interpolation} gives $o_P(i)\ne0$.
\end{proof}

\subsection{Binary second-degree obstructions in general DPOs}
\begin{defn}[Cone with a positivity gap]\label{def:gap-cone}
Let $P_0$ be a DPO positive cone containing $\Q_{\ge0}$, and set
\[
P=\{0\}\cup(1+P_0).
\]
All orders and monoid rings in this subsection are defined by this $P$.
\end{defn}

\begin{prop}\label{prop:gap-dpo}
The set $P$ is a DPO positive cone. For $a=2$, $p=1$, $q=3/2$, and $u=6$, the regions
\[
Y_1=\{v:1\le_Pv,\ 2v\le_P6\},\qquad
Y_2=\{v:3/2\le_Pv,\ 2v\le_P6\}
\]
have least elements $1,3/2$, respectively, and
\[
Y_1\cap Y_2=\{5/2,3\},
\]
The two points in the intersection are incomparable.
\end{prop}
\begin{proof}
For $r,s\in P_0$,
\[
(1+r)+(1+s)=1+(1+r+s),\qquad
(1+r)(1+s)=1+(r+s+rs),
\]
so $P$ is closed under addition and multiplication. If $1+r=-(1+s)$, then $-2=r+s\in P_0$, contradicting $2\in P_0$ and pointedness. Every $x=r-s$, with $r,s\in P_0$, equals $(1+r)-(1+s)$, so $P-P=L$.
Since $6-2\cdot1=4\in P$ and $6-2\cdot(3/2)=3\in P$, the two pieces have least elements $1,3/2$, respectively.
If $v\in Y_1\cap Y_2$, then $v\ne3/2$ because $1/2\notin P$. Hence $v-3/2\in1+P_0$, that is, $v\ge_{P_0}5/2$. If $v\ne3$, then $6-2v\in1+P_0$; since $\tfrac12P_0\subseteq P_0$, we obtain $v\le_{P_0}5/2$, and pointedness gives $v=5/2$.
Conversely, the quantities
\[
\begin{array}{c|ccc}
v&v-1&v-3/2&6-2v\\\hline
5/2&3/2&1&1\\
3&2&3/2&0
\end{array}
\]
all belong to $P$, so the intersection consists exactly of the two listed points. Neither their difference $1/2$ nor $-1/2$ belongs to $P$, so they are incomparable.
\end{proof}

\begin{thm}[A binary second-degree obstruction]\label{thm:gap-tor}
For the cone in Def~\ref{def:gap-cone},
\[
\Tor_2^{\Z[P]}\bigl(\Z[P]/(X^1,X^{3/2}),\Z[P]_2\bigr)_6\cong\Z.
\]
\end{thm}
\begin{proof}
Let $K_i=\ord(Y_i)$. Both pieces are upward closed in $Y_1\cup Y_2$, so each chain lies in a piece containing its least vertex. Hence
\[
\ord(Y_1\cup Y_2)=K_1\cup K_2,\qquad
K_1\cap K_2=\ord(\{5/2,3\}).
\]
Each $K_i$ is a cone, and their intersection consists of two discrete points. The Mayer--Vietoris segment
\[
0\longrightarrow\widetilde H_1(K_1\cup K_2)
\longrightarrow\widetilde H_0(K_1\cap K_2)
\longrightarrow0
\]
gives $\widetilde H_1(K_1\cup K_2)\cong\Z$. Taking $V=(1+P)\cup(3/2+P)$, $a=2$, and $u=6$ in Thm~\ref{thm:tor-order} gives the asserted $\Tor_2$.
\end{proof}

\begin{cor}\label{cor:gap-not-lattice}
The cone $P$ in Def~\ref{def:gap-cone} is not a lattice. For $Y_1,Y_2$ in Prop~\ref{prop:gap-dpo}, let $N$ be the nerve of the cover $\{\ord(Y_1),\ord(Y_2)\}$. Then $N\cong\Delta^1$, and
\[
\widetilde H_1(|N|;\Z)=0,\qquad
\widetilde H_1(Y_1\cup Y_2;\Z)\cong\Z.
\]
\end{cor}
\begin{proof}
If $P$ were a lattice, Cor~\ref{cor:binary-high-vanish} would force the $\Tor_2$ in Thm~\ref{thm:gap-tor} to vanish, a contradiction. Since $Y_1,Y_2,Y_1\cap Y_2$ are nonempty, $N$ is a full edge. The Mayer--Vietoris calculation in the proof of Thm~\ref{thm:gap-tor}, together with Prop~\ref{prop:mccord}, gives the latter homology group.
\end{proof}

\section{Realization and vanishing of higher obstructions}
\subsection{Realization of finite simplicial complexes}
\begin{thm}[Positive-sum realization of finite complexes]\label{thm:realization}
Let $P$ be a lattice containing finite families of pairwise reduced nonzero positive elements of arbitrarily large size. Every finite simplicial complex $K$ containing the empty face is isomorphic to some $\Delta_u(\mathbf p)$.
\end{thm}
\begin{proof}
If $K=\{\varnothing\}$, take $p>_P0$ and $u=0$, giving $\Delta_0(p)=\{\varnothing\}$.
Otherwise, label its vertex set by $[m]$, and let $\mathcal N$ be the set of minimal nonfaces. Choose pairwise reduced nonzero positive elements $e_0,e_N$ for $N\in\mathcal N$, and set
\[
p_i=e_0+\sum_{N\ni i}e_N,\qquad
u=me_0+\sum_{N\in\mathcal N}(|N|-1)e_N.
\]
For $I\subseteq[m]$,
\[
u-\sum_{i\in I}p_i
=(m-|I|)e_0+\sum_{N\in\mathcal N}(|N|-1-|I\cap N|)e_N.
\]
Cor~\ref{cor:disjoint-coordinates} gives
\[
I\in\Delta_u\iff\forall N\in\mathcal N,\ |I\cap N|<|N|
\iff\forall N\in\mathcal N,\ N\not\subseteq I\iff I\in K.
\]
For $\mathcal N=\varnothing$, the same formula gives the full simplex.
\end{proof}

\begin{cor}[Torsion in integral homology]\label{cor:torsion}
Under the hypotheses of Thm~\ref{thm:realization}, there are a finite family of positive elements $\mathbf p$ and $u\in P$ such that
\[
H_2K(\mathbf p)_u\cong\Z/2,\qquad
H_2\bigl(K(\mathbf p)\otimes_{\Z}\Q\bigr)_u=0.
\]
In particular, graded Koszul homology may contain nonzero torsion.
\end{cor}
\begin{proof}
Take a finite triangulation $K$ of $\mathbb{RP}^2$. Its cellular chain complex is
\[
0\longrightarrow\Z\xrightarrow{\,2\,}\Z\xrightarrow{\,0\,}\Z\longrightarrow0,
\]
so $\widetilde H_1(K;\Z)=\Z/2$. By Thm~\ref{thm:realization} and Thm~\ref{thm:delta}, choose $\Delta_u(\mathbf p)\cong K$ to obtain
\[
H_2K(\mathbf p)_u\cong\widetilde H_1(K;\Z)\cong\Z/2.
\]
Since $\Q$ is flat and $(\Z/2)\otimes\Q=0$, changing to rational coefficients annihilates this class.
\end{proof}

\subsection{Spheres in multiplicative homology}
\begin{defn}[Common-positive-part configuration]\label{def:sphere-config}
Let $P$ be a lattice. A common-positive-part configuration consists of $a,p_1,\ldots,p_m,b$ such that $m\ge2$, $a,b>_P0$, the elements $p_i>_P0$ are pairwise reduced, and
\[
b\le_Pap_i\qquad(1\le i\le m).
\]
Write $s=\sum_i p_i$, $u=as-b$, and $J=(X^{p_1},\ldots,X^{p_m})$.
\end{defn}

\begin{thm}[Realization of multiplicative spheres]\label{thm:sphere}
Every common-positive-part configuration satisfies $u\in P$ and
\[
\Xi_{a,u}(\mathbf p)=\partial\Delta^{m-1},\qquad
|\ord(Y_{a,u}(\cup_i(p_i+P)))|\simeq S^{m-2}.
\]
Consequently,
\[
\Tor_{m-1}^R(R/J,R_a)_u\cong\Z.
\]
\end{thm}
\begin{proof}
Since $ap_1-b\in P$, we have $u=(ap_1-b)+\sum_{i>1}ap_i\in P$. Reducedness and Cor~\ref{cor:disjoint-coordinates} give $b_I=\sum_{i\in I}p_i$. For $I\subsetneq[m]$, choose $j\notin I$. Then
\[
u-ab_I=(ap_j-b)+\sum_{k\notin I,\ k\ne j}ap_k\in P;
\qquad u-ab_{[m]}=-b\notin P.
\]
Thus $\Xi_{a,u}=\partial\Delta^{m-1}$, whose realization is $S^{m-2}$. Thm~\ref{thm:nerve} gives the homotopy equivalence of order complexes, and Thm~\ref{thm:xi} gives
\[
\Tor_{m-1}^R(R/J,R_a)_u
\cong\widetilde H_{m-2}(S^{m-2};\Z)\cong\Z.
\]
\end{proof}

\begin{cor}[Constructing a common positive part from a reduced family]\label{cor:sphere-polynomial}
Let $P$ be a lattice. Given $m\ge2$ pairwise reduced nonzero positive elements $p_i$, set
\[
b=\prod_{i=1}^mp_i,\qquad
a=\sum_{i=1}^m\prod_{j\ne i}p_j
\]
These form a configuration as in Def~\ref{def:sphere-config}.
\end{cor}
\begin{proof}
Products of nonzero field elements are nonzero, and pointedness ensures that finite sums of nonzero positive elements are nonzero. Hence $a,b>_P0$. For each $i$,
\[
ap_i-b=\sum_{k\ne i}p_i\prod_{j\ne k}p_j\in P,
\]
so $b\le_Pap_i$, verifying all configuration conditions.
\end{proof}

\subsection{Positive translations and vanishing of homology classes}
\begin{thm}[A specified positive translation annihilates the sphere class]\label{thm:death}
For a common-positive-part configuration, $\Xi_{a,u+b}(\mathbf p)$ is a full simplex. Therefore
\[
X^b:\Tor_{m-1}^R(R/J,R_a)_u
\longrightarrow\Tor_{m-1}^R(R/J,R_a)_{u+b}
\]
is the zero map.
\end{thm}
\begin{proof}
For every $I\subseteq[m]$,
\[
u+b-ab_I=a\sum_{j\notin I}p_j\in P.
\]
Thus $\Xi_{a,u+b}$ is a full simplex and has zero reduced homology. By Thm~\ref{thm:xi}, the target group is zero, and $X^b$ corresponds to the map induced by $\Xi_{a,u}\hookrightarrow\Xi_{a,u+b}$.
\end{proof}

\begin{defn}[Flat dimension]\label{def:flat-dimension}
The flat dimension $\fd_RM$ of a left $R$-module $M$ is the least length of a flat resolution; it is $\infty$ if no finite flat resolution exists. In particular, $\fd_RR_a$ always refers to the left action in Def~\ref{def:twist}.
\end{defn}

\begin{cor}[Unbounded obstructions for a fixed multiplier]\label{cor:infinite-fd}
If $a>_P0$ is fixed and common-positive-part configurations of length $m$ exist for every $m\ge2$, then $\fd_RR_a=\infty$.
\end{cor}
\begin{proof}
For each $m\ge2$, Thm~\ref{thm:sphere} gives a module $N_m$ with
\[
\Tor_{m-1}^R(N_m,R_a)\ne0.
\]
If a flat resolution $0\to F_d\to\cdots\to F_0\to R_a\to0$ of length $d$ existed, splitting it into short exact sequences and using $\Tor_{j\ge1}^R(N,F_k)=0$ would successively give
\[
\Tor_n^R(N,R_a)\cong\Tor_{n-d}^R(N,F_d)=0\quad(n>d).
\]
Taking $m>d+1$ gives a contradiction.
\end{proof}

\begin{rem}\label{rem:two-models}
Thm~\ref{thm:realization} realizes the complex $\Delta$ defined by positive sums, whereas Thm~\ref{thm:sphere} realizes $\Xi$, defined by multiplication and suprema. They correspond to different chain complexes; the former theorem does not imply that every finite complex can be realized as the latter.
\end{rem}

\section{Bao--Zhang and Hahn lattices}
\subsection{Structure of the Bao--Zhang lattice}
This section takes the lattice construction of \cite[\S5]{BY} as input. Fix all its choices and denote the resulting positive cone by $P_*$. We extract only the properties needed for the subsequent homology calculations.

\begin{defn}[Notation for the height construction]\label{def:height-data}
Let $L=\overline{\Q(T)}$ with $T\ne\varnothing$. Choose $t\in T$ and $i\in L$ with $i^2=-1$. Use the valuation height $d=-v$ of \cite[\S4]{BY}, with $d(0)=-\infty$, and choose compatible rational powers $t^s$. Then
\[
d(xy)=d(x)+d(y),\qquad d(t^s)=s,\qquad
d(x+y)\le\max\{d(x),d(y)\}.
\]
Write $L_{<h}=\{x\in L:d(x)<h\}$ and
\[
H_s=\{ct^s+x:c\in\Q_{>0},\ x\in L_{<s}\},\qquad
p=(t+i)/2,\quad q=(t-i)/2.
\]
\end{defn}

\begin{prop}[Data supplied by the construction]\label{prop:construction-data}
The cone $P_*$ has the following properties:
\begin{enumerate}
\item It is a lattice containing $\Q_{\ge0}$, and $1,p,q$ are pairwise reduced.
\item Every nonrational nonzero positive element belongs to some $H_s$ with $s\ge1$. Thus
\[
\{a\in P_*\setminus\{0\}:a^{-1}\in P_*\}=\Q_{>0};
\]
\item At every stage with finite cutoff height $A$, the lattice-ordered subspace $U$ satisfies $U\subseteq\{d\le A\}$. At the successor stage, a height $h>A$ is chosen, and the coordinate direct sum is extended by
\[
e_0=t^h,\qquad e_j=t^h+x_j,
\]
where $(x_j)$ represents a basis of $L_{<h}/U$. These new coordinate directions are pairwise reduced and remain so in the final cone.
\item If a nonrational positive element $a\in H_s$ already belongs to the lattice-ordered subspace at that stage, then $ae_j-1\in P_*$ for every new direction above.
\end{enumerate}
\end{prop}
\begin{proof}
For (1), see \byref{13}{Lem~5.1(ii)} and \byref{16}{Thm~5.4}; for (2), see \byref{16}{Thm~5.4} and \byref{19}{Remark~6.3}; for (3), see \byref{14}{Prop~5.2} and \byref{8}{Thm~3.4}.
For (4), let $S$ be the height monoid at that stage. From $a\in H_s$, $a\in U$, and the invariants of the Bao--Zhang construction, we have $s\in S\setminus\{0\}$. The successor height $h$ gives $S'=S+\N_0h$, and the successor cone contains every $H_r$ with $r\in S'$ and $r>h$; see \byref{14}{Prop~5.2}. By \byref{11}{Lem~4.2},
\[
ae_j\in H_sH_h=H_{s+h},\qquad d(1)=0<s+h,
\]
so $ae_j-1\in H_{s+h}\subseteq P_*$.
\end{proof}

\begin{cor}[New reduced families of arbitrary finite size]\label{cor:new-directions}
Let $U$ be a stage subspace in the Bao--Zhang construction with finite cutoff height. For every $m\ge1$, the next successor stage contains $e_1,\ldots,e_m\in P_*\setminus U$ that are pairwise reduced in $P_*$. If $a\in(P_*\cap U)\setminus\Q$, these elements can be chosen with
\[
ae_j\ge_{P_*}1\qquad(1\le j\le m).
\]
\end{cor}
\begin{proof}
Let the current cutoff be $A$ and the successor height be $h>A$. Choose distinct $r_1,\ldots,r_N\in\Q\cap(A,h)$. For $c_j\in\Q$ not all zero, the valuation inequality and the distinct heights give
\[
d\!\left(\sum_jc_jt^{r_j}\right)=\max_{c_j\ne0}r_j>A.
\]
Thus these elements are linearly independent in $L_{<h}/U$, so this quotient is infinite dimensional. Prop~\ref{prop:construction-data}(3) gives new reduced families of arbitrary finite size; if $a\in U$, (4) gives $ae_j\ge_{P_*}1$.
Every $a\in L$ belongs to a stage subspace with finite cutoff height; see \byref{16}{Thm~5.4}.
\end{proof}

\subsection{Multiplicative homology and flat dimension}
\begin{thm}[An explicit first-degree class]\label{thm:explicit-p}
Set $R_*=\Z[P_*]$, $J=(X^p,X^q)$, and $u=p^2+pq-1$. Then $u\in P_*$, and
\[
\Tor_1^{R_*}(R_*/J,(R_*)_p)_u\cong\Z.
\]
The multiplication map
\[
X^1:\Tor_1^{R_*}(R_*/J,(R_*)_p)_u
\longrightarrow\Tor_1^{R_*}(R_*/J,(R_*)_p)_{u+1}
\]
is zero.
\end{thm}
\begin{proof}
The initial cone contains $H_2$; see \byref{13}{Lem~5.1}. Since $d(i)=d(1)=0$,
\[
p^2-1=(t^2+2it-5)/4\in H_2,\qquad
pq-1=(t^2-3)/4\in H_2.
\]
Thus $1\le_{P_*}p^2,pq$, while $p,q$ are pairwise reduced. In Def~\ref{def:sphere-config}, take $a=p$, $(p_1,p_2)=(p,q)$, and $b=1$, giving $u=a(p_1+p_2)-b$. Thm~\ref{thm:sphere} and Thm~\ref{thm:death} give the isomorphism and annihilation by $X^1$, respectively.
\end{proof}

\begin{thm}[Flat dimensions for $P_*$]\label{thm:original-fd}
For every $a\in P_*\setminus\{0\}$,
\[
\fd_{R_*}(R_*)_a=
\begin{cases}
0,&a\in\Q_{>0},\\
\infty,&a\notin\Q_{>0}.
\end{cases}
\]
More precisely, in the nonrational case, for each $n\ge1$ there exist pairwise reduced elements $e_1,\ldots,e_{n+1}$ such that
\[
\Tor_n^{R_*}\!\left(R_* /(X^{e_1},\ldots,X^{e_{n+1}}),(R_*)_a\right)_{a\sum_j e_j-1}
\cong\Z.
\]
\end{thm}
\begin{proof}
If $a\in\Q_{>0}$, then $a^{-1}\in P_*$, and Thm~\ref{thm:flat-mult} gives $\fd_{R_*}(R_*)_a=0$.
If $a\notin\Q_{>0}$, fix $n\ge1$. By Cor~\ref{cor:new-directions}, choose $n+1$ pairwise reduced nonzero positive elements $e_j$ satisfying
\[
ae_j-1\in P_*\quad(1\le j\le n+1).
\]
Apply Thm~\ref{thm:sphere} with $b=1$ to obtain
\[
\Tor_n^{R_*}\!\left(R_* /(X^{e_1},\ldots,X^{e_{n+1}}),(R_*)_a\right)_{a\sum_j e_j-1}
\cong\widetilde H_{n-1}(S^{n-1};\Z)\cong\Z.
\]
Since $n$ is arbitrary, Cor~\ref{cor:infinite-fd} gives infinite flat dimension.
\end{proof}

\begin{cor}\label{cor:fixed-p-all}
Fix $p=(t+i)/2$. For every $n\ge1$, there are pairwise reduced nonzero positive elements $e_1,\ldots,e_{n+1}\in P_*$. Set
\[
J_n=(X^{e_1},\ldots,X^{e_{n+1}}),\qquad
u_n=p\sum_{j=1}^{n+1}e_j-1\in P_*,
\]
Then
\[
\Tor_n^{R_*}(R_*/J_n,(R_*)_p)_{u_n}\cong\Z,
\]
and the map induced by $X^1$ from this component to degree $u_n+1$ is zero.
\end{cor}
\begin{proof}
Since $d(p)=1$, we have $p\notin\Q$. Taking $a=p$ in Thm~\ref{thm:original-fd} gives nonzero classes in every positive degree. These configurations all have $b=1$, so Thm~\ref{thm:death} gives vanishing after the upper bound increases by $1$.
\end{proof}

\subsection{Construction of the Hahn lattice}
\begin{defn}[Coefficientwise positive cone]\label{def:hahn}
Take a transcendental element $s$, the countable field $K=\overline{\Q(s)}$, and the cone $P_0$ given by the Bao--Zhang construction on $K$. Set
\[
H=K((z^{\Q})),\qquad
P_H=\left\{\sum_{r\in\Q}c_rz^r:c_r\in P_0\text{ for all }r\right\}.
\]
The nonzero support of each series in $H$ is required to be well ordered in the usual order on $\Q$.
\end{defn}

\begin{thm}[Hahn lattice]\label{thm:hahn}
The cone $P_H$ is a lattice on $H$, with lattice operations defined coefficientwise. There exists a field isomorphism $H\cong\C$. For any such isomorphism $\varphi:H\to\C$, the image $\varphi(P_H)$ is a lattice on $\C$.
\end{thm}
\begin{proof}
The sum of well-ordered supports is well ordered, and only finitely many decompositions of a fixed exponent contribute to a product coefficient, so Hahn multiplication is well defined. For $f,g\in P_H$, every product coefficient
\[
[z^r](fg)=\sum_{s+t=r}([z^s]f)([z^t]g)
\]
belongs to $P_0$, so $P_HP_H\subseteq P_H$. Additive closure and pointedness follow coefficientwise from $P_0$.
For $f=\sum_rc_rz^r$, set $f^\pm=\sum_rc_r^\pm z^r$. Since $\supp f^\pm\subseteq\supp f$, these are still Hahn series, and
\[
f=f^+-f^-,\qquad f^+,f^-\in P_H.
\]
Hence $P_H-P_H=H$. For any $g=\sum_rd_rz^r\in H$, the two series
\[
(f\vee g)_r=c_r\vee d_r,\qquad(f\wedge g)_r=c_r\wedge d_r
\]
have support contained in $\supp f\cup\supp g$, which is well ordered. The universal properties of coefficientwise suprema and infima give lattice operations on the whole field.
By \cite[Theorem 1]{MacLane}, $H$ is algebraically closed because the coefficient field $K$ is algebraically closed and the exponent group $\Q$ is divisible. Moreover,
\[
2^{\aleph_0}\le |H|\le |K|^{|\Q|}=2^{\aleph_0},
\]
where the lower bound comes from ordinary power series with coefficients in $\{0,1\}$. Thus $\operatorname{trdeg}_{\Q}H=2^{\aleph_0}=\operatorname{trdeg}_{\Q}\C$. A bijection between transcendence bases extends to their algebraic closures, giving $H\cong\C$.
\end{proof}

\begin{defn}[Unobstructed multipliers]\label{def:units}
For a DPO positive cone $P$, write
\[
\mathcal U(P)=\{a\in P\setminus\{0\}:a^{-1}\in P\}.
\]
By Thm~\ref{thm:flat-mult}, these are exactly the positive multipliers inducing flat changes of coefficients. If nonempty, this set is a group under field multiplication.
\end{defn}

\begin{prop}[Unobstructed multipliers of the Hahn cone]\label{prop:hahn-units}
We have
\[
\mathcal U(P_H)=\{cz^r:c\in\Q_{>0},\ r\in\Q\}.
\]
In particular, $z$ is transcendental over $K$, and $1,z\in\mathcal U(P_H)$ satisfy $1\wedge z=0$.
\end{prop}
\begin{proof}
If $f,g\in P_H\setminus\{0\}$ and $fg=1$, every product of nonzero positive coefficients is nonzero, and pointedness prevents cancellation in a finite sum of positive terms. Hence
\[
\supp f+\supp g=\supp(fg)=\{0\}.
\]
Fix $s\in\supp g$. Every $r\in\supp f$ satisfies $r+s=0$, so $\supp f$ is a singleton; likewise for $\supp g$. Thus
\[
f=cz^r,\qquad g=c^{-1}z^{-r},\qquad c,c^{-1}\in P_0.
\]
By Prop~\ref{prop:construction-data}(2), $c\in\Q_{>0}$. Conversely, these monomials and their inverses are coefficientwise positive. The element $z$ is transcendental over $K$ because monomials with distinct integer exponents are linearly independent. The identity $1\wedge z=0$ follows from coefficientwise infima.
\end{proof}

\begin{prop}[Unbounded multiplicative obstructions in the new cone]\label{prop:hahn-infinite}
Let $R_H=\Z[P_H]$. Then
\[
\fd_{R_H}(R_H)_{1+z}=\infty.
\]
\end{prop}
\begin{proof}
Fix $n\ge1$ and take the $2n$ pairwise reduced positive elements
\[
e_j=z^{3j},\quad f_j=z^{3j+1}\quad(0\le j<n),\qquad
u_n=\sum_{j=0}^{n-1}(z^{3j}+z^{3j+1}+z^{3j+2}).
\]
For a subset $I$ of these directions, let $\alpha_j,\beta_j\in\{0,1\}$ record whether $e_j,f_j$ belong to $I$. Reducedness gives $b_I=\sum_j(\alpha_je_j+\beta_jf_j)$, so
\[
u_n-(1+z)b_I
=\sum_j\bigl((1-\alpha_j)z^{3j}
+(1-\alpha_j-\beta_j)z^{3j+1}+(1-\beta_j)z^{3j+2}\bigr).
\]
The exponents are distinct, and $P_0\cap\Q=\Q_{\ge0}$. Thus the difference belongs to $P_H$ if and only if $\alpha_j+\beta_j\le1$ for each $j$. Hence
\[
\Xi_{1+z,u_n}=(S^0)^{*n},\qquad
|\Xi_{1+z,u_n}|\cong S^{n-1}.
\]
By Thm~\ref{thm:xi}, for $J_n=(X^{e_j},X^{f_j}:0\le j<n)$,
\[
\Tor_n^{R_H}(R_H/J_n,(R_H)_{1+z})_{u_n}\cong\Z.
\]
Since $n$ is unbounded, dimension shifting along a flat resolution gives $\fd_{R_H}(R_H)_{1+z}=\infty$.
\end{proof}

\subsection{Inequivalence under field automorphisms}
\begin{defn}[Equivalence under field automorphisms]\label{def:equivalence}
For two DPO positive cones $P,Q$ on $\C$, define
\[
P\sim Q\iff\exists\sigma\in\operatorname{Aut}(\C):\ \sigma(P)=Q.
\]
\end{defn}

\begin{thm}[Two inequivalent lattices]\label{thm:inequivalent}
For any field isomorphism $\varphi:H\to\C$, set $Q=\varphi(P_H)$. Then $Q$ is inequivalent to the cone $P_*$ given by the Bao--Zhang construction on $\C$.
\end{thm}
\begin{proof}
If $\sigma\in\operatorname{Aut}(\C)$ satisfied $\sigma(P_*)=Q$, then $\psi=\varphi^{-1}\circ\sigma:\C\to H$ would satisfy
\[
\psi(P_*)=P_H,\qquad
\psi(\mathcal U(P_*))=\mathcal U(P_H),\qquad
\psi(x\wedge y)=\psi(x)\wedge\psi(y).
\]
But $\mathcal U(P_*)=\Q_{>0}$, where any two elements $r,s$ satisfy $r\wedge s=\min\{r,s\}>0$, whereas $1,z\in\mathcal U(P_H)$ satisfy $1\wedge z=0$, a contradiction.
\end{proof}

\begin{cor}\label{cor:not-all-extreme}
There exist a lattice $P$ on $\C$ and $a\in P\setminus\Q$ such that
\[
\fd_{\Z[P]}(\Z[P])_a=0.
\]
\end{cor}
\begin{proof}
Choose a field isomorphism $\varphi:H\to\C$, and set $P=\varphi(P_H)$ and $a=\varphi(z)$. By Prop~\ref{prop:hahn-units}, $z$ is transcendental over $\Q$ and $z,z^{-1}\in P_H$, so $a\notin\Q$ and $a,a^{-1}\in P$. Thm~\ref{thm:flat-mult} gives the asserted flatness.
\end{proof}

\begin{ques}\label{question:finite-positive}
Do there exist a lattice $P$ on $\C$ and $a>_P0$ such that
\[
0<\fd_{\Z[P]}\Z[P]_a<\infty?
\]
The constructions above give flat modules and modules of infinite flat dimension, but do not settle this question.
\end{ques}

\begin{rem}[Scope of the existence problem]\label{rem:existence}
The reducedness obstructions developed here give criteria for a candidate positive cone to be a lattice. An independent proof of the existence of a lattice on $\C$ using these criteria would still require an extension theorem for compatible elimination of obstructions. The existence input in this section remains the construction in \cite{BY}.
\end{rem}

\section*{Acknowledgments}
I am especially grateful to Professor Zhang YueHui of Shanghai Jiao Tong University. He introduced me to the Birkhoff--Pierce problem and to the broad mathematical landscape behind it. Together we wrote \cite{BY}. Throughout this work, he has patiently explored every question with me, guided me back when I lost my way, and continually encouraged me to keep exploring.

His guidance has extended far beyond individual problems and proofs. He has helped me cultivate mathematical taste and understand research from a wider perspective. Like his renowned course \emph{The Sky of Mathematics}, he has offered me glimpses of the beauty of many fields of mathematics and opened a larger world to me. I count myself fortunate to have met him.

\section*{AI Collaboration Statement}
This paper grew out of a sustained collaboration between the author and GPT6Astra. The author primarily supplied the overall framework, the ideas behind each proposed result, the directions for mathematical experimentation, and rough outlines of possible proofs. Following the author's directions and requirements, GPT6Astra primarily checked arguments, supplied missing details, constructed counterexamples, and searched the literature. The work developed through successive layers of mutual influence; neither contribution can realistically be disentangled from the other. As stated in the opening note, this manuscript is a rough draft written throughout by the model under the author's natural-language guidance. The author bears full responsibility for the correctness of its contents.

\end{document}